\documentclass[runningheads]{svjour2b}
\smartqed  
\usepackage{graphicx}
\newcommand \R {{I\!\!R}} 
\newcommand \Z {{\mathrm Z}} 
\begin{document}

\title{On uniqueness of semi-wavefronts }

\subtitle{Diekmann-Kaper theory of a nonlinear convolution equation re-visited}


\author{Maitere Aguerrea, Carlos Gomez   and Sergei Trofimchuk
}


\institute{M. Aguerrea  \at
              Facultad de Ciencias B\'asicas, Universidad Cat\' olica del Maule, Casilla 617, Talca  \ 
              \email{maguerrea@ucm.cl}           
           \and
           C. Gomez and S. Trofimchuk  \at
          Instituto de Matem\'atica y Fisica, Universidad de Talca, Casilla 747,
Talca, Chile \  E-mails: cgomez@inst-mat.utalca.cl,
trofimch@inst-mat.utalca.cl
}

\date{\notused}

\maketitle

\begin{abstract}
Motivated by the uniqueness problem for monostable semi-wave\,-fronts, we propose a revised version of the  Diekmann and Kaper theory of a nonlinear convolution equation.   
Our version of the Diekmann-Kaper theory allows 1) to consider new types of models which include nonlocal KPP type equations (with either symmetric or anisotropic dispersal), nonlocal lattice equations and delayed reaction-diffusion equations;   2) to incorporate the critical case (which corresponds to the slowest wavefronts) into the consideration;  3)   to  weaken or to remove various restrictions on  kernels and nonlinearities. The results are compared with those of Schumacher (J. Reine Angew. Math.  316:  54-70, 1980),  Carr and Chmaj (Proc. Amer. Math. Soc. 132:  2433-2439, 2004), and other more recent studies.   
{\bf Mathematics Subject ClassiÞcation (2010)} 45J10, 45J05, 35K57, 35C07 
\end{abstract}

\section{Introduction}
\label{intro}
The main goal of this paper is to develop a version of the fundamental Diekmann and Kaper theory \cite{D1,D2,dk} (the DK theory for short) of a nonlinear convolution equation  for the scalar integral equation
\begin{equation}\label{17} \hspace{-7mm}
\varphi(t) = \int_Xd\mu(\tau)\int_\R K(s,\tau)g(\varphi(t-s),\tau) ds, \quad t \in \R, 
\end{equation}
in the case of  monostable nonlinearity  $g$. 
Throughout the paper $(X,\mu)$ will denote a measure space with  finite measure $\mu$, $K(s,\tau) \geq 0$  will be
integrable on $\R\times X$ with $ \ \int_{\R}K(s,\tau)ds>0, \ \tau \in X,$ while measurable $g: \R_+ \times X \to \R_+,$ $g(0,\tau) \equiv 0, $ will be continuous 
in $\varphi$ for every fixed $\tau \in X$.  When $X$ is just a single point (i.e. $\# X =1$), 
equation (\ref{17}) coincides with  the nonlinear convolution equation from  \cite{dk}.

In a biological context, $\varphi$ is the size of an adult population, so we are interested in non-negative solutions of (\ref{17}).  Following the terminology of \cite{GK},  we call a bounded continuous non-constant solution $\varphi: \R \to \R_+$ semi-wavefront   if either
$\varphi(-\infty) = 0$ or $\varphi(+\infty) = 0$. We will always assume  $\varphi$ to satisfy
$\varphi(-\infty) =0$, since  the other case 
can be easily transformed to this one via the change of variables
$\zeta(t)= \varphi(-t)$,  with equation  (\ref{17}) assuming the form  
$$
\zeta(t) = \int_Xd\mu(\tau)\int_\R K_1(s,\tau)g(\zeta(t-s),\tau) ds, \quad K_1(s,\tau):=K(-s,\tau).$$  We would like to emphasize that the nonlinearity $g$ and semi-wavefronts are generally non-monotone \cite{FT} (nevertheless,  typically semi-wavefronts are strictly increasing  in some vicinity  of $-\infty$ \cite{atv,fhw,TT}).  The non-monotonicity of waves complicates their analysis. For instance, the wave uniqueness  is easier to establish within a subclass of monotone solutions \cite{CDM,GT,wlr2}. 

Actually the  `largely open uniqueness question' \cite{CG}   is central  in our research where we follow the scheme elaborated in  \cite{dk}. This means that after assuming the existence of a semi-wavefront  to (\ref{17}),  we study its asymptotic behavior at infinity trying then to demonstrate the  wave uniqueness (modulo  translation).  Similarly to other authors,  we work mostly with the first positive eigenvalue $\lambda_l$ of the linearization of  (\ref{17}) at zero. As a consequence,  our analysis excludes from the consideration so called "pushed" fronts \cite{ES,GK,SA} associated to the second positive eigenvalue $\lambda_r$. 
Analogously to  \cite{dk},  the existence of semi-wavefronts to (\ref{17})  is not investigated here.

There are various motivations to study the above equation, mainly from the theory of traveling waves for nonlinear 
models  (e.g. reaction-diffusion equations with delayed response \cite{atv,GT,tz,TAT,wlr2}, equations with non-local dispersal \cite{BFR,CC,Co,CDM,KPP,KS}, lattice systems \cite{CG,fwz,GW,MaZ}).  Only a few of these models take
the simplest form  with $\# X =1$ of (\ref{17}).  Therefore our first goal is to show that the basic framework
of \cite{dk} can be extended to include much broader class of  convolution type equations than it was initially intended. Here is a simple step to create such a general  direct extension of results in \cite{dk}.  It would be interesting to consider  further generalizations of (\ref{17}) in order to include more applications (for example, equations with distributed delays considered in  \cite{fwz,fz},  see also \cite{GH,KS,wlr2}).  However, we do not pursue this direction in our current work.  After all, ours is not the first attempt to expand the DK theory.   Schumacher   has mentioned, while studying equation 
$$
c\varphi'(t) = g(\varphi, \mu_c * g(\varphi)), 
$$
the impossibility of transforming it into the form  to which the DK theory could be applied \cite[p.54]{KS}.  Instead, Schumacher has developed an approach which is based on guidelines of the DK theory and, at the same time,  which is technically rather different from that in \cite{dk}.  In particular,  in order to extend the DK uniqueness theorem,  Schumacher has used a comparison method for differential inequalities combined with Nagumo-point argument.   In this respect,  his work \cite{KS} is very close to the recent  contributions \cite{CG,Co,CDM,GW}.   

Similarly to \cite{KS},  the present studies also follow the mainstream of the DK ideology.  In difference with \cite{KS} and trying to apply our  results to delayed equations (where in general the comparison argument does not work), we preserve the original idea of the DK theory in the proof of   uniqueness.  Now, from the technical point of view our approach to equation (\ref{17}) differs  from the methods used by  Diekmann and Kaper, Schumacher and Carr and Chmaj \cite{CC} in many key points.  Even though the logical sequence of results here  basically is the same  as in \cite{dk}, our proofs are essentially different. In  particular,   we do not  use  the Titchmarsh theory of Fourier integrals \cite{dk,fwz} nor we use the  Ikehara Tauberian theorem \cite{CC,CDM,wlr2} in order to obtain asymptotic expansions of solutions (a necessary key component of each uniqueness proof). We have found more convenient  for our purpose 
the use of  a suitable $L^2-$variant  of  the bootstrap argument (as it was suggested by Mallet-Paret in \cite[p. 9-10]{FA}).

As a consequence of the DK strategy,  we  also present a non-existence result and describe properties of the kernel $K$  which is proved to satisfy  exponential convergence estimates (Mollison's condition \cite{CDM}).   Here  the fulfillment of the Mollison's condition means that the characteristic function
$$
 \chi(z) := 1 - \int_Xg'(0,\tau)d\mu(\tau)\int_ {\R}
K(s,\tau)e^{-zs}ds
$$ is well defined for all $z$ from some maximal
non-degenerate interval (which can be open, closed, half-closed, finite or infinite).  
One of the key results of the theory says that, under rather mild assumptions on $g, K$ the presence of a semi-wavefront $\varphi, \varphi(-\infty)=0,$ guarantees the existence of a minimal positive zero $\lambda_l$ to $\chi(z)$.    The spreading properties of  some integro-differential equations with `fat-tailed' kernels were recently considered by Garnier \cite{Ga}. 

Next, as it is known the DK and Schumacher uniqueness theorems do not apply to the critical fronts (when $\chi(\lambda_l)= \chi'(\lambda_l)=0$).   As an example, let us consider 
 the nonlocal KPP equation 
 \begin{equation}\label{ccd}
 u_t = J * u - u+ g(u), \ x \in \R, \  g(0) = g(1) =0, \ f >0 \ {\rm on \ } (0,1)
\end{equation}
proposed in \cite{KPP}.
Here  continuous birth function $f$  is supposed  to be differentiable  at $0$, with  $g(s) = g'(0)s + O(s^{1+\alpha}), \ s \to 0+$,  for some $\alpha >0$, and to satisfy the KPP condition \cite{KPP} $f'(s) \leq f'(0)$,  $  s \in (0,1)$.  Measurable kernel $J\geq 0, \int J ds =1$ is allowed to be asymmetric  and non compactly supported. This agrees with  the initial idea of Kolmogorov, Petrovsky and Piskunov \cite{KPP} who 
interpreted $J(x)dx$ as the probability that an individual passes a distance between $x$ and $x+dx$.  
It is easy to see that the DK theory does not apply to (\ref{ccd}).  Under the above mentioned assumptions, Schumacher \cite[Example 2]{KS} has proved uniqueness of all {\it non-critical} wavefronts for (\ref{ccd}). 
Later on, Carr and Chmaj  \cite{CC}  achieved  an important extension of the  DK theory for the special case of equation (\ref{ccd}).  By assuming several additional conditions in \cite{CC} that  $J$ must be even, compactly supported and 
 \begin{equation}\label{lcs}
 |g(s)-g(t)|\leq g'(0)|t-s|, \quad s,t \geq 0, 
 \end{equation}
they showed that 
 the minimal wavefront $\varphi(x+c_0t)$ to (\ref{ccd}) satisfying  $ 0 \leq \varphi(s) \leq 1, \ s \in \R, $   is unique up to translation.  
Carr and Chmaj's work  has motivated the second goal of our research:  to get an improvement of the DK theory  that includes the critical semi-wavefronts.  Theorem   \ref{mth1} below gives such an extension for general model (\ref{17}).  In the particular case of equation (\ref{ccd}) our result (stated as Theorem \ref{th5}) establishes the uniqueness of {\it critical} wavefronts under the same assumptions on $J, f$ as in \cite{KS}.  
See  Section  6.1 for more details, further discussion and references. 

The necessity of the subtangential Lipcshitz condition (\ref{lcs})  \cite{CC,dk,fwz,tz}  could be considered as a weak point of the DK uniqueness theorem, cf. \cite{atv,CG,CDM,GW,GK,KS}.   For instance, as it was established recently by Coville, D\'avila and Mart\'inez \cite{CDM},  neither (\ref{lcs}) nor $g'(s)\leq g'(0),$ $  s \in (0,1),$ is necessary  to prove the uniqueness of  non-stationary monotone traveling fronts to (\ref{ccd}). Instead of that, it was supposed in \cite{CDM} that generally asymmetric  $J \in C^1(\R)$  is compactly supported with $J(a) >0, $ $J(b) >0$ for some $a<0<b$,  while $g \in C^1(\R)$ has to satisfy
$g'(0)g'(1) <0, $ $  g(s) \leq g'(0)s, \ s \geq 0$, and $g \in C^{1, \alpha}$ near $0$.  The proof in \cite{CDM} follows ideas of \cite{Co} and is mainly based on the sliding methods proposed by Berestycki and Nirenberg \cite{BN} (see \cite{Co,CDM} for a comprehensive state-of-art overview about (\ref{ccd}) and \cite{Chen,MaZ} for the further references).    The above discussion explains our third goal in this paper: to weaken various convergence and smoothness conditions of the DK theory, and especially condition (\ref{lcs}).   It is worthwhile to note that a similar task was also considered in \cite{KS}.   The related improvements  can be found  in Theorems \ref{mth1} and \ref{mth2}. In the latter theorem,  we remove condition (\ref{lcs}) by  assuming a little more smoothness for $g$ and exploiting the absence of zeros for $\chi(z)$ in the vertical strip  $\lambda_l < \Re z < \lambda_r$ (see Lemma \ref{tg}).  
Incidentally,  Theorems \ref{mth2} justifies the following principle  for monostable equations: 
"fast positive semi-wavefronts are unique (modulo translation)".  In the last section, we  apply this principle to reaction-diffusion equations with delayed Mackey-Glass type nonlinearities.

The main results of this paper are stated as Theorems \ref{mth1},  
\ref{mth2} below.  We apply them to nonlocal integro-differential equations (Section 6.1),  nonlocal lattice systems (Section 6.2),  nonlocal (Section 6.3) and local (Section 6.4) reaction-diffusion equations with discrete delays. 
In Theorem \ref{thm:expnv}, we give  a short proof of  the necessity of 
the  Mollison's condition  for the existence of semi-wavefronts. Theorem \ref{gg}  provides a non-existence result.

\vspace{-2mm}
\section{Mollison's condition}
\vspace{-2mm} 

In this section, we consider somewhat more general  equation 
\begin{equation} \label{tweabnv}
\varphi(t) = \int_Xd\mu(\tau)\int_\R K(s,\tau)g(\varphi(t-s),t-s,\tau) ds,
\end{equation}
where measurable $g: \R\times\R \times X \to \R_+$ is continuous 
in the first two variables for every fixed $\tau \in X$. We suppose additionally that, for some 
measurable $p(\tau) \geq 0$  and $ \delta >0, \ \bar s  \leq 0$, it holds 
\begin{equation}\label{gco}
 g(v,s, \tau) \geq p(\tau)v,  \ v \in (0,\delta),  \ s \leq \bar s, \ \tau \in X.  
\end{equation}
First, we present a simple proof of the necessity of the following Mollison's condition (cf. \cite{CDM}) 
for the existence of the semi-wavefronts:

\vspace{-3mm}

\begin{eqnarray}\label{kkk}\hspace{-7mm}
\int_\R \int_XK(s,\tau)p(\tau)d\mu(\tau)e^{-sz}ds \ {\rm is \ finite \ for
\ some } \ z \in \R\setminus\{0\}.
\end{eqnarray}
\vspace{-3mm}
\begin{theorem}\label{thm:expnv}
Let continuous $\varphi:\R\rightarrow[0,+\infty)$ satisfy (\ref{tweabnv}) and suppose that $\varphi(-\infty)=0$ and 
$\varphi(t)\not\equiv 0, \ t \leq t'$ for each fixed $t'$.  If (\ref{gco}) holds and 
\begin{equation} \label{condinv} \int_X\int_\R K(s,\tau)p(\tau)dsd\mu(\tau) \in (1,\infty), 
\end{equation}
then $\int_{-\infty}^0\varphi(s)e^{-s\bar x}ds $ and $\int_\R \int_XK(s,\tau)p(\tau)d\mu(\tau)e^{-s\bar x}ds$ are convergent for an appropriate   $\bar x>0$. Furthermore, supp\,$K \cap (\R_+\times X) \not=\emptyset$. 
\end{theorem}
\begin{remark} Looking for   heteroclinic solutions of the simple  logistic equation $x' = - \beta x + x(1+\beta-x)$ with $\beta >0$,  we obtain an example of  (\ref{17}) where  supp\,$K \cap (\R_-\times X) =\emptyset$ under conditions of  the above theorem.  
\end{remark}
\begin{proof} Since the support of $K$ generally is  unbounded, we will truncate $K$  by 
choosing integer  $N$ such that 
$$\kappa:=\int_X\int_{-N}^N K(s,\tau)p(\tau)dsd\mu(\tau) >1, \  {\rm and \ } 0\leq  \varphi(t) <\delta, \ t < \bar s - N.$$ Integrating equation (\ref{tweabnv}) between $t'$ and $t< \bar s -N$, we find that
\begin{eqnarray*}\hspace{-9mm}
\int_{t'}^{t}\varphi(v)dv\geq &&
 \int_Xd\mu(\tau)\int_{-N}^NK(s,\tau)\int^t_{t'}g(\varphi(v-s),v-s,\tau) dvds\\ \geq
&& \int_Xp(\tau)d\mu(\tau)\int_{-N}^NK(s,\tau)\int^t_{t'}\varphi(v-s)dvds
 \\
=&& \int_Xp(\tau)d\mu(\tau)\int_{-N}^NK(s,\tau)(\int^{t'}_{t'-s}+ \int^{t}_{t'}+\int^{t-s}_{t})\varphi(v)dvds,\end{eqnarray*}
from which 
$$
\int_{t'}^{t}\varphi(v)dv\leq \frac{2\delta \int_X\int_{-N}^N |s|K(s,\tau)p(\tau)dsd\mu(\tau)}{\int_X\int_{-N}^N K(s,\tau)p(\tau)dsd\mu(\tau)-1}, \quad t' < t < \bar s -N. 
$$
Hence,  increasing function 
\begin{equation} \label{pov}
\psi(t) = \int_{-\infty}^t \varphi(s)ds
\end{equation}
is well defined for all $t \in \R$ and  
$$
\psi(t)
\geq \int_Xp(\tau)d\mu(\tau)\int_{-N}^NK(s,\tau)\psi(t-s)ds \geq 
\kappa \psi(t-N),  \quad t < \bar s -N. 
$$
Consider $h(t)=\psi(t)e^{-\gamma t}$ where ${\kappa}=e^{\gamma N}$, cf. \cite{CC}.  For all $t< \bar s -N$ we have $$h(t-N)=\psi(t-N)e^{-\gamma(t-N)}\leq\frac{1}{\kappa}\psi(t)e^{-\gamma t}e^{\gamma N}=h(t)$$ and $\gamma={N^{}}\ln{{\kappa}}>0$. Hence $\sup_{t \leq 0}{h(t)}<\infty$ and $\psi(t)=O(e^{\gamma t})$, $t\rightarrow-\infty$. After taking $\bar x \in (0,\gamma)$ and integrating by parts, we obtain
$$
\int_{-\infty}^t\varphi(s)e^{-\bar xs}ds = \psi(t)e^{-\bar xt}+ \bar x \int_{-\infty}^t\psi(s)e^{-\bar xs}ds
$$ 
that proves the first statement of the theorem. 
Finally, 
$$
e^{-\bar x t}\psi(t)
= \int_Xd\mu(\tau)\int_\R e^{-\bar x s}K(s,\tau)e^{-\bar x (t-s)}\psi_1(t-s,\tau)ds, 
$$ 
where $\psi_1(u, \tau) := \int_{-\infty}^ug(\varphi(s),s,\tau)ds \geq p(\tau) \int_{-\infty}^u\varphi(s)ds, \ u \leq \bar s -N.$ The latter
yields 
$$
\int_{-\infty}^{\bar s-N}e^{-\bar x v}\psi(v)dv
= \int_Xd\mu(\tau)\int_{\R}e^{-\bar x s}K(s,\tau)\int_{-\infty}^{\bar s-N }e^{-\bar x (v-s)}\psi_1(v-s,\tau)dvds\geq 
$$
$$
 \int_Xp(\tau)d\mu(\tau)\int_{-\infty}^{0}e^{-\bar x s}K(s,\tau)ds\int_{-\infty}^{\bar s -N}e^{-\bar x v}\psi(v)dv,  $$
\begin{equation} \label{polk}
\mathcal{K}_-(\bar x):=\int_Xp(\tau)d\mu(\tau)\int_{-\infty}^{0}e^{-\bar x s}K(s,\tau)ds\leq 1, \  ({\rm note \ that  \ } \psi(s) >0, \ s \in \R),
\end{equation}
so that 
$$
\mathcal{K}_-(0)=\int_Xp(\tau)d\mu(\tau)\int_{-\infty}^{0}K(s,\tau)ds\leq 1 < \int_Xp(\tau)d\mu(\tau)\int_{\R}K(s,\tau)ds, 
$$
which completes the proof of  the theorem. 
\hfill {}
\end{proof}
\begin{remark} Suppose that $|g(\varphi(s),s,\tau)| \leq C$ where $C$ does not depend on $s,\tau$.   
Then 
$$
|\varphi(t+h)- \varphi(t)| \leq C\int_\R|K_a(s +h)- K_a(s)|ds, $$
where $K_a(s) := \int_XK(s, \tau)d\mu(\tau) \in L_1(\R)$. Since the translation is continuous in $L_1(\R)$ \cite[Example 5.4]{en}, we find that $\varphi(t)$ is uniformly continuous on $\R$. 
It is easy to see that the convergence of the integral $\int_{-\infty}^0\varphi(s)ds < \infty$ combined with the uniform continuity  of $\varphi$  gives   $\varphi(-\infty) =0$. In this way,  $\int_{-\infty}^0\varphi(s)ds < \infty$  implies that $\int_{-\infty}^0e^{-xs}\varphi(s)ds < \infty$  for small positive $x$.  
\end{remark} 
\begin{remark}  It is easy to see that the global  non-negativity of $g$ is not necessary  
in the case of $K$ having bounded support (uniformly in $\tau \in X$). 
\end{remark}

Now, let $\varphi, K, g, \bar x $ be as in Theorem \ref{thm:expnv}. Set 
$$\displaystyle \Phi(z)=\int_{\R}e^{-z s}\varphi(s)ds,\  \mathcal{K}(z)= \int_\R \int_XK(s,\tau)p(\tau)d\mu(\tau)e^{-sz}ds,$$ and denote the maximal open vertical strips of convergence for these two integrals as $\sigma_\phi < \Re z < \gamma_\phi$ and $\sigma_K < \Re z < \gamma_K$, respectively. Evidently, $\sigma_\phi, \sigma_K \leq 0$ and $\gamma_\phi, \gamma_K \geq \bar x >0$. 
 Since $\varphi, K$ are both non-negative, by \cite[Theorem 5b, p. 58]{widder},   $\gamma_\phi, \ \gamma_K, \sigma_\phi, \ \sigma_K$ are singular points  of  $\Phi(z), \mathcal{K}(z)$ (whenever they  are finite). 
 A simple inspection of the proof of Theorem \ref{thm:expnv} suggests  the following 
 \begin{lemma} \label{fn}  Assume $\varphi, g, K$ are as in Theorem \ref{thm:expnv}.  Then
 $\sigma_K \leq \sigma_\phi < \gamma_\phi \leq \gamma_K$.  Furthermore, 
 $\mathcal{K}(\gamma_\phi )$ is always a finite number. 
 \end{lemma}
 \begin{proof} 
 For all $z \in (0, \gamma_\phi), \ t \leq 0,$ we have
 $$
 \psi(t) = \int_{-\infty}^t(\varphi(s) e^{-zs})e^{zs}ds \leq e^{zt}\int_{-\infty}^0\varphi(s) e^{-zs}ds,
 $$
 so that $\int_{-\infty}^0 \psi(s)e^{-z's}ds < \infty$ for each $z' \in (0,\gamma_\phi)$ and, due to (\ref{polk}),  we get 
 $$
 \mathcal{K}_-(z):=\int_Xp(\tau)d\mu(\tau)\int_{-\infty}^{0}e^{-zs}K(s,\tau)ds\leq 1 $$ 
for all $z \in (0, \gamma_\phi)$.  Hence, using the Beppo Levi monotone convergence theorem, we obtain that $\mathcal{K}_-(\gamma_\phi)\leq 1$. As a  consequence, $\mathcal{K}(\gamma_\phi)$ is finite and  $\gamma_K \geq \gamma_\phi$.  
 \end{proof}
 \begin{corollary}\label{cfi} Assume that 
 $$
 \lim_{z \to \gamma_K-} \int_\R \int_XK(s,\tau)p(\tau)d\mu(\tau)e^{-sz}ds =+ \infty. 
 $$
 Then $\gamma_\phi$ is a finite number and $\gamma_\phi < \gamma_K$. 
 \end{corollary}

\vspace{-2mm}
\section{Abscissas of convergence}
\vspace{-2mm} 
In this section, we investigate the abscissas of convergence  for the bilateral Laplace transforms of  $K$ and bounded non-negative  $\varphi$ satisfying $\varphi(-\infty)=0$, $\varphi(t)\not\equiv 0, \ t \leq t',$ for each fixed $t'$, and solving our main equation (\ref{17}). 
Now we are supposing that the continuous $g(\cdot,\tau): \R_+ \to \R_+$ is differentiable at $0$ with  $g'(0+,\tau) >0$ for each fixed $\tau$.  Then 
the non-negative functions
$$
\lambda^+_\delta(\tau):= \sup_{u \in (0,\delta)}\frac{g(u,\tau)}{u},\ \lambda^-_\delta(\tau):= \inf_{u \in (0,\delta)}\frac{g(u,\tau)}{u}, \quad \delta >0, \  \tau \in X,
$$
are well defined, measurable, monotone in $\delta$ and pointwise converging:
$$
\lim_{\delta \to 0+} \lambda^\pm_\delta(\tau) = g'(0+,\tau). 
$$  
The {\it characteristic} function $\chi$ associated  with the variational equation along the trivial 
steady state of (\ref{17}) is defined by  
$$
\chi(z):= 1-\int_\R \int_XK(s,\tau)g'(0+,\tau)d\mu(\tau)e^{-sz}ds.
$$ 
It is supposed to be negative at $t =0$:  $\chi(0) <0$. 
Since condition (\ref{gco}) is obviously satisfied with $p(\tau)= \lambda^-_\delta(\tau)$ and 
$$
\lim_{\delta \to 0+} \int_\R \int_XK(s,\tau)\lambda^-_\delta(\tau)d\mu(\tau)ds=\int_\R \int_XK(s,\tau)g'(0+,\tau)d\mu(\tau)ds >1
$$ 
by the monotone convergence theorem,  all  results of Section 2 hold true for equation (\ref{17}).
Furthermore, we have the following
\begin{theorem} \label{gg}  Assume $\chi(0)<0$.  Let $\varphi:\R\rightarrow[0,+\infty)$ be a semi-wavefront  to  equation (\ref{17}). If $\varphi(-\infty)=0$  and $\varphi(t)\not\equiv 0, \ t \leq t'$ for each fixed $t'$, then $\chi(z)$  has a zero on $(0, \gamma_\phi]\subset (0, \gamma_K]\subset \R\cup \{+\infty\}$.  
\end{theorem}
\begin{remark}  1) If  $\varphi(+\infty)=0$  then a similar statement can be proved. Namely, in such a case  $\chi(z)$  has a zero on $[\sigma_K,0)$.  2) It should be noted that Theorem \ref{gg} also provides a non-existence result: 
if $\chi(x) <0$ for all $x \in (0,\gamma_K]$ then equation (\ref{17}) does not have any semi-wavefront vanishing at $-\infty$. 
\end{remark}
\begin{proof} 
For real positive $z \in (0, \gamma_\phi)$ we consider the integrals 
$$
\Phi(z)=\int\limits_{\R}e^{-z s}\varphi(s)ds,  \mathcal{G}(z,\tau):=\int\limits_{\R}e^{-z s}g(\varphi(s),\tau)ds, 
\mathcal{K}(z,\tau):=\int\limits_{\R}e^{-z s}K(s,\tau)ds. 
$$
Since $\varphi$ is non-negative and bounded, and since $g'(0+,\tau)>0$ exists, the convergence 
of $\mathcal{G}(z,\tau)$ (for positive $z$) is equivalent to the convergence of $\Phi(z)$. 
Applying  the bilateral Laplace transform to equation (\ref{17}),  we obtain that
\begin{equation}\label{chtos}
\Phi(z)= \int_X\mathcal{K}(z,\tau)\mathcal{G}(z,\tau)d\mu(\tau).
\end{equation}
Obviously, $\mathcal{K}, \mathcal{G}, {\Phi}$
are positive at each real point of the convergence. 

Let us prove that $\chi(z)$  has a zero on $(0, \gamma_\phi]$.   First, we suppose that $\Phi(\gamma_\phi)= \lim_{z \to \gamma_\phi-}\Phi(z) = \infty$.  In such a  case, we claim that 
$$
\lim_{z \to \gamma_\phi-}\frac{\mathcal{G}(z,\tau)}{\Phi(z)} = g'(0,\tau). 
$$
Indeed, let $T_\delta$ be the rightmost non-positive number such that $\varphi(s) \leq \delta$ for $s \leq T_\delta$.  Then 
$$\lambda^-_\delta \int_{-\infty}^{T_\delta}e^{-z s}\varphi(s)ds \leq \int_{-\infty}^{T_\delta}e^{-z s}g(\varphi(s),\tau)ds \leq \lambda^+_\delta \int_{-\infty}^{T_\delta}e^{-z s}\varphi(s)ds,
$$
$$
 \int^{+\infty}_{T_\delta}e^{-z s}(g(\varphi(s),\tau)+\varphi(s))ds\leq \frac{\sup_{s \in \R}(g(\varphi(s),\tau)+\varphi(s))}{z}e^{-\gamma_\phi T_\delta}. 
$$
As a consequence, for each positive $\delta >0$, 
$$
\lambda^-_\delta \leq \liminf_{z \to \gamma_\phi-}\frac{\mathcal{G}(z,\tau)}{\Phi(z)} \leq \limsup_{z \to \gamma_\phi-}\frac{\mathcal{G}(z,\tau)}{\Phi(z)} \leq \lambda^+_\delta, 
$$
that proves our claim.

Now, by using the Fatou lemma as  $z \to \gamma_\phi-$ in  
$$
\int_X\mathcal{K}(z,\tau)\frac{\mathcal{G}(z,\tau)}{
\Phi(z)} d\mu(\tau)=1, 
$$
we obtain 
$$
1-\chi(\gamma_\phi) = \int_X\mathcal{K}(\gamma_\phi,\tau)g'(0,\tau) d\mu(\tau) \leq 1.  
$$
Therefore
 $\chi(\gamma_\phi)\geq 0,$ and since  $\chi(0) <0$ we get the required assertion.  

Hence, we may suppose that  $\Phi(\gamma_\phi)= \lim_{z \to \gamma_\phi-}\Phi(z)>0$ is finite. Since $\varphi(t)\not\equiv 0, \ t \leq t'$ for each fixed $t'$, in such a case $\gamma_\phi < \infty$.  Due to  Lemma  \ref{fn}, the value $\mathcal{K}(\gamma_\phi)$ is also finite. Set  $$\zeta(t): = \varphi(t)e^{-\gamma t}, \  K_1(s,\tau):= e^{-\gamma s}K(s,\tau), \ {\rm  where\ } \ \gamma:= \gamma_\phi.$$ Then, for $t < T_\delta -N$, we have from (\ref{17}) that
$
\int_{-\infty}^t\zeta(v)dv =
$
$$
\int_{-\infty}^t\varphi(v)e^{-\gamma v}dv \geq \int_Xd\mu(\tau)\int_{-N}^N K_1(s,\tau)\int_{-\infty}^tg(\varphi(v-s),\tau) e^{-\gamma(v-s)}dvds \geq 
$$
$$
 \int_Xd\mu(\tau)\int_{-N}^N\lambda^-_\delta(\tau) K_1(s,\tau)\int_{-\infty}^t\zeta(v-s)dvds \geq 
$$
$$
( \int_Xd\mu(\tau)\int_{-N}^N\lambda^-_\delta(\tau) K_1(s,\tau)ds) \int_{-\infty}^{t-N}\zeta(v)dv.  
$$
Suppose now on the contrary that the characteristic equation 
$$
\chi(z):= 1- \int_\R \int_XK(s,\tau)g'(0+,\tau)d\mu(\tau)e^{-sz}ds=0 
$$
has not real roots on $[0,\gamma_\phi]$. Then $\chi(0)<0$ implies $\chi(\gamma) <0$.  As a consequence, in virtue of the monotone convergence theorem, 
$$
 \lim_{\delta \to 0+, N \to +\infty}\int_Xd\mu(\tau)\int_{-N}^N\lambda^-_\delta(\tau) K_1(s,\tau)ds= 1-\chi(\gamma) > 1.  
$$
Hence, for some appropriate $\delta, N >0$,  increasing function 
$
\xi(t) = \int_{-\infty}^t \zeta(s)ds
$
satisfies
$
\xi(t)
 \geq \kappa_\delta \xi(t-N), $ $   t < T_\delta -N
$
with $\kappa_\delta >1$. Arguing now as in the proof of Theorem \ref{thm:expnv} below (\ref{pov}) 
we conclude that the integral 
$\int_{-\infty}^t\zeta(s)e^{-zs}$ converges for all small positive $z$,  contradicting to the definition of $\gamma_\phi$. \hfill {}
\end{proof}
\begin{remark} It is clear that $\chi(z)$ is concave on $(\sigma_K, \gamma_K)$, where $\chi''(z) <0$. 
Since $\chi(0)$ is negative,  $\chi$ can have at most two real zeros, and they must be  of the same sign.   We will denote them (if they exist) by $\lambda_l \leq \lambda_r$. Under assumption of the existence of a semi-wavefront $\varphi$ vanishing at $-\infty$, $\chi$ has at least one positive root $\lambda_l$.  Finally, it is clear that $\chi$ is analytical in the vertical strip $\Re z \in (0, \gamma_K)$. 
\end{remark}

\noindent {\bf Notation} At this stage, it is convenient to introduce the following notation: 
$$
 \lambda_{rK}=\left\{%
\begin{array}{ll}
    \lambda_r, & \hbox{if} \  \lambda_r \   \hbox{exists},    \\
     \gamma_K, & \hbox{otherwise}.
\end{array}%
\right.
$$

\begin{lemma} \label{tg} Equation $\chi(z)=0$ does not have roots in the open strip $\Sigma:= \Re z \in (\lambda_l, \lambda_{rK})$.  Furthermore, the only possible zeros on the boundary $\Sigma$ are $\lambda_l, \lambda_r$. 
\end{lemma}
\vspace{-5mm}
\begin{proof}
Observe that if $\chi(z_0)=0$ for some $z_0 \in \Sigma$, then 
$\chi(\Re z_0) > 0$ since $\chi$ is concave, $\chi(\lambda_l)=0$ and 
$\Re z_0 \in (\lambda_l, \min\{\lambda_r, \gamma_K\})$.  On the other hand, $1=$
$$
 |\int_\R \int_XK(s,\tau)g'(0+,\tau)d\mu(\tau)e^{-sz_0}ds| \leq \int_\R \int_XK(s,\tau)g'(0+,\tau)d\mu(\tau)e^{-s\Re z_0}ds
$$
and therefore $\chi(\Re z_0) \leq 0$, a contradiction. Now, if $\chi(\lambda_l + i\omega) =0$ for some 
$\omega \not= 0$ then similarly
$$
1=  \chi(\lambda_l + i\omega) = |\chi(\lambda_l + i\omega)| \leq \chi(\lambda_l)=1,  
$$
so that   
$$
\int_\R \int_XK(s,\tau)g'(0+,\tau)d\mu(\tau)e^{-s\lambda_l}(1- \cos \omega s)ds=0. 
$$
Thus $K(s,\tau)(1- \cos \omega s) =0$ for almost all $\tau \in X$, so that $K(s,\tau) =0$ a.e. on $X \times \R$, 
a contradiction. 
  \hfill {}
\end{proof}
\vspace{-2mm}
\section{A bootstrap argument}
\vspace{-2mm} 
The main purpose of this section is to prove several auxiliary statements 
needed in the studies of the asymptotic behavior of solutions $\varphi(t)$ at $t=-\infty$. 
Usually proofs of the uniqueness are based on the derivation of appropriate 
asymptotic formulas with one or two leading terms (at $t=-\infty$ as in \cite{CC,dk,fwz,wlr2} or at $t=+\infty$ as in \cite{GT}).    
 Our approach is based on an asymptotic integration routine often used in the theory of functional  differential equations, e.g. see \cite{Hale},  \cite[Proposition 7.1]{FA} or 
\cite{FT2}.   Thus we use neither the Titchmarsh theory of Fourier integrals \cite{Tit} nor the powerful Ikehara Tauberian theorem \cite{CC,dk}.  First we will apply our methods to get an asymptotic formula for the integral $\psi(t) := \int_{-\infty}^t\varphi(s)ds$. Since $\psi \in C^1(\R)$ is strictly increasing and positive, this function is somewhat easier 
to treat than the  solution $\varphi(t)$. 

Everywhere in the sequel, we assume all conditions of Section 3 on $\varphi, K, g, \chi$. In particular, $\chi(0)<0$.  
We also will use  the following hypotheses ({\bf{SB}}), ({\bf{EC}}$_\rho$): 
\begin{description}
\item[ {\bf{ (SB)}}] $\gamma_\phi < \gamma_K$ and, for some measurable  $C(\tau) >0$ and $\alpha, \sigma \in (0,1]$,
$$
|g'(0,\tau)-\frac{g(u,\tau)}{u}| \leq C(\tau)u^{\alpha}, \ u \in (0,\sigma), $$
\begin{equation}\label{ze}
\zeta(x) := \int_{X\times \R}C(\tau)K(s,\tau)e^{-sx}dsd\mu <  +\infty, \ x \in (0, \gamma_K). 
\end{equation}
\end{description}
\begin{description}
\item[ {\bf{ (EC$_\rho$)}}] For every $x \in (0, \rho), \ \rho \leq \gamma_\phi,$ there exists some positive $C_x$ such that 
\begin{equation} \label{UEC}
0 \leq \varphi(t) \leq C_xe^{xt},  \ t \leq 0.  
\end{equation} 
\end{description}
There are several situations when  {\bf{ (EC$_\rho$)}} can be easily checked:
\begin{lemma} \label{eq} Condition {\bf{ (EC$_\rho$)}} is satisfied in either of the following two cases: 
\begin{itemize}
\item[(i)]  $\varphi \in C^1(\R)$ and the integral $ \ \int_{\R}e^{-xs}\varphi'(s)ds$ converges absolutely for all $x  \in (0, \rho)$;
\item[(ii)] (cf. \cite{dk}) $\rho < \gamma_\phi$ and  there exist  measurable 
$d_1, d_2, \  d_1d_2 \in L^1(X), $ such that  
$$
0\leq K(s,\tau) \leq d_1(\tau)e^{\rho s}, \ s \in \R,\  \tau \in X,  
$$
\begin{equation}\label{d2}
|g(u,\tau)| \leq d_2(\tau)u, \ u \geq 0. 
\end{equation}
\end{itemize}
\end{lemma} 
\vspace{-3mm}
\begin{proof} {\it (i)} For each  $x \in (0, \rho)$  we have that  
$$
\varphi(t) = \int_{-\infty}^t\varphi'(s)ds = \int_{-\infty}^te^{xs}\varphi'(s)e^{-xs}ds \leq e^{xt}  
\int_{\R}e^{-xs}|\varphi'(s)|ds=: C_x e^{xt}. 
$$
{\it (ii)} Since $\rho < \gamma_\phi$, the integral  $\int_{\R}e^{-xs}\varphi(s)ds$ converges for all $x \in (0,\rho]$.  If  $x \in (0, \rho], \ t \leq 0$,  then
$$
\varphi(t)e^{-x t} \leq \varphi(t)e^{-\rho t} = \int_Xd\mu(\tau)\int_{\R}K(s,\tau)e^{-\rho s}e^{-\rho(t-s)}g(\varphi(t-s),\tau) ds \leq 
$$
$$
C:=\int_Xd_1(\tau)d_2(\tau)d\mu(\tau)\int_\R e^{-\rho s }\varphi(s)ds,  \ t \in \R.   
$$
\end{proof}
The following simple proposition will be used several times in the sequel: 
\begin{lemma} \label{uc} Assume that $h(s)e^{-sx} \in L^1(\R)$ for all $x \in [a,b]$. Then 
$$
H(x,y): = \int_{\R}h(s)e^{-sx-isy}ds, \ y \in \R,
$$
is uniformly (with respect to  $y \in \R$) continuous on $[a,b]$. 
\end{lemma}
\begin{proof} Take an arbitrary $\varepsilon  >0$ and 
let $N>0$ be such that
$$
\int_{\R\setminus [-N,N]}|h(s)|e^{-sx}ds < 0.25 \varepsilon ,  \  x \in [a,b]. 
$$ 
Since $e^t$ is uniformly continuous on compact sets, there exists $\delta >0$ such that 
$|x_1-x_2| \leq \delta, \ s \in [-N,N]$ implies $|e^{-x_1s}- e^{-x_2s}| < 0.5\varepsilon/|h|_1$. But then 
$$
|H(x_1,y)-H(x_2,y)| \leq 0.5 \varepsilon + \int_{-N}^N|h(s)||e^{-x_1s} -e^{-x_2s}|ds < \varepsilon, \ y \in \R.  
$$
\end{proof}
\begin{corollary} \label{inn} With $h$ as in Lemma \ref{uc}, we have that 
$
\lim_{y \to \infty}H(x,y) =0 
$
uniformly on $x \in [a,b]$. 
\end{corollary}
\begin{proof} Due to Lemma \ref{uc},  for each $\varepsilon>0$  there exists 
a finite sequence $a:=x_0< x_1<x_2< \dots < x_m=:b$ possessing the following property:  
for each $x$ there is  $x_j$ such that 
$|H(x_j,y)-H(x,y)| < 0.5\epsilon$ uniformly on $y$. Now, by the Riemann-Lebesgue lemma, 
$\lim_{y \to \infty}H(x_j,y) =0$ for every $j$. Therefore, for all $j$ and some $M>0$, we have that 
$|H(x_j,y)| <   0.5\epsilon$ if $|y|\geq M$. This implies that 
$$
|H(x,y)| \leq |H(x_j,y)-H(x,y)| +|H(x_j,y)| < \epsilon, \quad |y| \geq M,\  x \in [a,b],  
$$
and the corollary is proved.  \hfill {}
\end{proof}
As we know, the property $\varphi(-\infty) =0$ implies the exponential decay $\psi(t) = O(e^{z t})$ at $-\infty$  for each $z \in (0,\gamma_\phi)$.  It is clear also that $\psi(t) = O(t)$ as $t \to +\infty$. Hence,  for each fixed $z \in (0,\gamma_\phi)$,  we can integrate equation (\ref{17}) twice,   to find that $\Psi(z):=\int_{\R}e^{-zv}\psi(v)dv $ satisfies
$$
\Psi(z)=  \int_Xd\mu(\tau)\int_\R K(s,\tau)e^{-zs}\int_{\R}e^{-z(v-s)}\int_{-\infty}^{v-s}g(\varphi(u),\tau)dudvds=
$$
$$
\int_Xd\mu(\tau)\int_\R K(s,\tau)e^{-zs}\int_{\R}e^{-zv}\int_{-\infty}^{v}g(\varphi(u),\tau)dudvds=
$$
$$
\left(\int_Xd\mu(\tau)\int_\R K(s,\tau)g'(0,\tau)e^{-zs}ds\right)\int_{\R}e^{-zv}\psi(v)dv+ \mathcal{R}(z), \
{\rm \ where}
$$
$$
\mathcal{R}(z):=\int_Xd\mu(\tau)\int_\R K(s,\tau)e^{-zs}ds\int_{\R}e^{-zv}\int_{-\infty}^{v}(g(\varphi(u),\tau)-g'(0,\tau)\varphi(u))dudv.
$$
Therefore 
$
\chi(z)\Psi(z)= \mathcal{R}(z). 
$
Set now
$$
\mathbf{G}(z,\tau):=\int_{\R}e^{-zv}G(v,\tau)dv, \quad 
G(v, \tau):=\int_{-\infty}^{v}(g(\varphi(u),\tau)-g'(0,\tau)\varphi(u))du.
$$
\begin{lemma} \label{afp} Assume (\ref{d2}), ({\bf{SB}}), {\bf{ (EC$_{2\epsilon}$)}} for some small  $2\epsilon \in (0,  \gamma_K - \gamma_\phi)$. Then given $a,b \in (0, \gamma_\phi + \alpha \epsilon)$ there exists  $\rho>0$ 
depending on $\varphi,a,b$ such that 
$$
|\mathbf{G}(z,\tau)| \leq \rho(\tau)/|z|: = \rho(C(\tau)+d_2(\tau)+ g'(0,\tau))/{|z|}, \quad \Re z \in [a,b] \subset (0, \gamma_\phi + \alpha \epsilon).
$$ 
\end{lemma}
\begin{proof} 
For $x:=\Re z \in (0, \gamma_\phi + \alpha \epsilon), \ v \leq 0$, we have
$$
e^{-x v}|G(v, \tau)|\leq e^{-xv}C(\tau) \int_{-\infty}^{v}(\varphi(u))^{1+\alpha}du \leq e^{-xv}C_\epsilon^\alpha C(\tau)\psi(v)e^{\alpha \epsilon v},$$ 
so that $e^{-x \cdot}|G(\cdot, \tau)| \in L^1(\R)\cap L^2(\R). $
After integrating by parts, we obtain 
$$
\int_{-N}^Ne^{-zv}G(v,\tau)dv= \frac{G(-N,\tau)e^{zN}- G(N,\tau)e^{-zN}}{z} + $$
$$+\frac{1}{z}\int_{-N}^Ne^{-zu}(g(\varphi(u),\tau)-g'(0,\tau)\varphi(u))du.
$$ 
This yields 
$$
|\int_{\R}e^{-zv}G(v,\tau)dv|= \frac{1}{|z|}|\int_{\R}e^{-zu}(g(\varphi(u),\tau)-g'(0,\tau)\varphi(u))du|\leq
$$ 
$$
\frac{1}{|z|}\left(C_\epsilon^\alpha C(\tau)\int\limits_{-\infty}^0e^{-(\Re z-\alpha\epsilon)u}\varphi(u)du + |\varphi|_\infty (g'(0,\tau) + d_2(\tau)) \int\limits^{+\infty}_0e^{-\Re zu}du \right). 
$$
 \end{proof}
 \begin{corollary} \label{kg} In addition, assume that $\int_{\R\times X}K(s,\tau)\rho(\tau) e^{-sx}d\mu ds$
 converges for all $x \in (0, \gamma_K)$. Then $\chi(\gamma_\phi) =0$ and, for  appropriate $\varepsilon_1>0, \ a,m  \in \R$, $k \in \{0,1\}$, and continuous $r \in L^2(\R)$,  it holds that
$$
\psi(t+m) = (a-t)^ke^{\gamma_\phi t} + e^{(\gamma_\phi+\varepsilon_1) t}r(t), \ t \in \R. 
$$
 \end{corollary}
It should be noted that depending on the geometric properties of $g$, the value of $\gamma_\phi$ can be minimal 
(the case of a pulled semi-wavefront \cite{ES,GK,SA}) or maximal (the case of a pushed semi-wavefront, ibid.)
positive zero of $\chi(z)$.  Observe that,  due to the monotonicity of $\psi$, we can also use here the Ikehara Tauberian  theorem \cite{CC}. However  it gives a slightly different result. 
\begin{proof}
Set $z: =x +i y$. For a fixed  $0< x <  \gamma_\phi+\alpha \epsilon$ we have 
$$
|\mathcal{R}(z)| = |\int_X\mathbf{G}(z,\tau)\int_\R K(s,\tau)e^{-zs}dsd\mu| \leq  
 \frac{1}{|z|}\int_X\rho(\tau)\int_\R K(s,\tau)e^{-xs}dsd\mu, 
$$
so that $\mathcal{R}(z)$ is regular in the strip $0< \Re z< \gamma_\phi+\alpha \epsilon$.   Thus  we can deduce from $
\Psi(z)= \mathcal{R}(z)/\chi(z)$ that $\gamma_\phi=\gamma_\psi$ (e.g. see \cite[Lemma 4.4]{dk}, the definition of $\gamma_\psi$ is similar to that of $\gamma_\phi$) must be a positive zero of $\chi(z)$ and  $\Psi(\gamma_\phi) =\infty$.  
 It is clear that  $\mathcal{R}(x+i\cdot)$  is also bounded and square integrable on $\R$ (for each fixed $x$).  Take now $\gamma' , \gamma ''$ such that $0 < \gamma' < \gamma_\phi < \gamma'' < \gamma_\phi + \alpha\epsilon$.  Then we may shift the path of integration in the inversion formula for the Laplace transform (e.g. see \cite[p. 10]{FA}) to obtain 
$$
\psi(t) = \frac{1}{2\pi i}\int_{\gamma' -i\infty}^{\gamma'
+i\infty}e^{zt}\Psi(z)dz
=- {\rm Res}_{z = \gamma_\phi} \frac{e^{zt}\mathcal{R}(z)}{\chi(z)} + \frac{e^{\gamma'' t}}{2\pi i}\left\{
\int_{-\infty}^{+\infty}e^{ist}a_1(s)ds \right\},$$ 
where the first term is different from $0$ and  $a_1(s) = \mathcal{R}(\gamma'' + is)/\chi(\gamma'' + is)$ is square integrable on $\R$.  Here we recall that, by Corollary  \ref{inn}, $\lim_{y \to \infty} \chi(x+i y) =1$ uniformly on $x \in [\gamma', \gamma'']$. 
Since $\chi''(x) >0, \ x \in (0, \gamma_K)$ ,  for some $a,m\in \R$  we get
$\psi(t+m) = (a-t)^k e^{\gamma_\phi t} + e^{\gamma'' t} r(t)$.  \hfill {}
\end{proof}
\begin{lemma} \label{nond} Assume all  conditions of Lemma \ref{afp} except $\gamma_\phi < \gamma_K$.  If 
$$
1- \chi_1(x_0):=\int_\R \int_XK(s,\tau)d_2(\tau)d\mu(\tau)e^{-sx_0}ds \leq 1,  
$$
for some $x_0 \in (0, \gamma_K)$, 
then  $\gamma_\phi$ coincides with the minimal  positive zero $\lambda_l$ of $\chi(z)$.
\end{lemma}
\begin{proof} Since $d_2(\tau)\geq g'(0,\tau)$, we obtain that $x_0 \in [\lambda_l, \lambda_{rK}] $  and $\lambda_l  < \gamma_K$.  \underline{Case I:}  $\gamma_\phi < \gamma_K$.  Then,  by
Corollary \ref{kg}, we have $\chi(\gamma_\phi)=0$ so that $\gamma_\phi \in \{\lambda_l, \lambda_r\}$. 
Suppose that $\gamma_\phi > \lambda_l$, this implies $x_0\leq  \gamma_\phi=\lambda_r$.   We have 
$$
\Psi(z)= 
\left(\int_Xd\mu(\tau)\int_\R K(s,\tau)d_2(\tau)e^{-zs}ds\right)\int_{\R}e^{-zv}\psi(v)dv+ \mathcal{R}_1(z), \
{\rm \ where}
$$
$$
\mathcal{R}_1(z):=\int_Xd\mu(\tau)\int_\R K(s,\tau)e^{-zs}ds\int_{\R}e^{-zv}\int_{-\infty}^{v}(g(\varphi(u),\tau)-d_2(\tau)\varphi(u))dudv, 
$$
or, in a shorter form, 
\begin{equation}\label{chi1}
\chi_1(z)\Psi(z)= \mathcal{R}_1(z). 
\end{equation}
It is clear that $x_0= \gamma_\phi = \lambda_r > \lambda_l$ implies immediately that $g'(0,\tau)= d_2(\tau)$ a.e. on $X$ and 
that $\chi_1(z)=\chi(z), \  \mathcal{R}(z)=\mathcal{R}_1(z).$ As we have seen in the proof of Corollary \ref{kg}, 
this guarantees that $\mathcal{R}_1(x_0)$ is  a finite number. Of course, $\mathcal{R}_1(x_0)$ is also well defined if $x_0 < \gamma_\phi$.   Now, it is clear that 
$\mathcal{R}_1(x_0) \leq 0$ because  of  $g(u,\tau) \leq d_2(\tau)u, \ u \geq 0$.  We claim that,  in fact,  $\mathcal{R}_1(x_0) < 0$. Indeed, otherwise  $g(u,\tau)=d_2(\tau)u, \ u \geq 0,$ for almost all $\tau \in X$  that yields $d_2(\tau)=g'(0,\tau)$ and $ \mathcal{R}_1(z) \equiv 0$ leading to a contradiction: $\Psi(z) \equiv 0$ and $\psi(t) \equiv 0$. 

Now, from $\mathcal{R}_1(x_0) < 0, \Psi(x_0) > 0, \chi_1(x_0) \geq 0$, we 
deduce that $\Psi$ must have a pole at $x_0=\gamma_\phi  < \gamma_K$. But then $\chi_1(\gamma_\phi)= \chi(\gamma_\phi)$ implies  $\chi_1(z)\equiv \chi(z)$, $\mathcal{R}(z)=\mathcal{R}_1(z)$.
Hence, $\lambda_l < \lambda_r=x_0 < \gamma_K$ and    $\gamma_\phi =x_0$ is a simple pole of $\Psi$. Therefore we can proceed as in the proof of Corollary \ref{kg} taking  $0 < \gamma' < \gamma_\phi = \lambda _r < \gamma'' < \gamma_\phi + \alpha\epsilon$ to obtain  
$$
\psi(t) = \frac{1}{2\pi i}\int_{\gamma' -i\infty}^{\gamma'
+i\infty}e^{zt}\Psi(z)dz
= -{\rm Res}_{z = \lambda_r} \frac{e^{zt}\mathcal{R}(z)}{\chi(z)} + e^{\gamma'' t}r_1(t)=
$$
$$
=  Ae^{ \gamma_\phi t} + e^{\gamma'' t}r_1(t),
\quad {\rm
where } \
A:=- \frac{\mathcal{R}( \lambda_r)}{\chi'(\lambda_r)}<0, \ r_1 \in L^2(\R), 
$$
contradicting to the positivity of $\psi$.  \\
\underline{Case II:}  $\gamma_\phi = \gamma_K$.  Since $x_0 < \gamma_K=\gamma_\phi$ and $\mathcal{R}_1(x_0) < 0$, we similarly deduce from  (\ref{chi1}) that $x_0$ is a singular point of $\Psi(z)$, a contradiction.
\end{proof}

\vspace{-2mm}
\section{The uniqueness theorems}
\vspace{-2mm} 

To prove our uniqueness results we will need  
more strong  property  of $\varphi$ than the merely convergence of $\int_{\R}e^{-zs}\varphi(s)ds$ for all $\Re z \in (0, \gamma_\phi)$ (even combined, as in Section 4, with {\bf {(EC$_{\epsilon}$)}} for some small $\epsilon >0$).  
This property, assumed everywhere in the sequel,  is  {\bf {(EC$_{\gamma_\phi}$)}}.  The nonlinearity $g$ is supposed to satisfy the hypothesis  {\bf {(SB)}}. 

The following assertion is crucial for  extension of the Diekmann-Kaper theory on the critical case $\chi(\lambda_l)= \chi'(\lambda_l)=0$.  
\begin{lemma} \label{trud}Suppose that, for some  $a, b >\delta >0$,  continuous  $v: \R \to [0,1)$ satisfies  $v(t)=1 + O(e^{a t}),\  t \to - \infty$, $v(t)= O(e^{-b t}), \ t \to +\infty,$ and 
$$
v(t) \leq \int_{\R}N(s)v(t-s)ds,
$$
where measurable $N(s) \geq 0, \ s \in \R,$ is such that
$$
\int_{\R}N(s)ds =1, \quad \int_{\R}sN(s)ds =0, \quad \int_{\R}N(s)e^{x s}ds < \infty, \ {\rm \  for  \  all \ } |x| \leq \delta. 
$$
Then $v(t) \equiv 0$.
\end{lemma}
\begin{proof} First we observe that, without restricting the generality, we may assume that 
$v \in C^2(\R)$ with the finite norm $|v|_{C^2}:=\sup_{s \in \R, j =0,1}|v^{(j)}(s)|$. Indeed, if we set
$$
w(t) := \int_t^{t+1}v(s)ds,    \ t \in \R, 
$$
then $w \in C^1(\R)$ has the same properties as $v$, $|w'(t)| < 1, \ t \in \R$, while $v(t) \equiv 0$ if and only if 
$w(t) \equiv 0$. For instance, if $v(t) \leq ce^{-bt}$ for $t \geq t_0, \ b >0$, then $w(t) \leq ce^{-bt}\int_0^1e^{-bs}ds\leq ce^{-bt}, \ t \geq t_0$.  Furthermore, $w'(t) = v(t+1)-v(t)$ behaves as  $O(e^{a t})$ at $ - \infty$ and as  
$O(e^{-b t})$ at $ +\infty$.  

Applying the same procedure to $w$ once more, we 
obtain the desired smoothness property of $v$ with $v'(t), v''(t)$ satisfying
\begin{equation}\label{oo}
v'(t), v''(t) =  O(e^{a t}),\  t \to - \infty, \quad v'(t), v''(t) = O(e^{-b t}), \ t \to +\infty. 
\end{equation}
In any case, the bilateral Laplace transform $V(z)$ 
of $v(t)$ is well defined in the vertical strip $-b< \Re z <0$. 

Set now 
$$
f(t):= \int_{\R}N(s)v(t-s)ds - v(t) \geq 0.
$$
It follows from this definition that $0 \leq f(t) \leq 1- v(t)$ and therefore 
 $f(t)= O(e^{at}), \ t \to -\infty$. Additionally, using (\ref{oo}), we obtain, for $j=0,1,2$ and some positive $C,C'>0$,  
$$
 \int_{\R}N(s)|v^{(j)}(t-s)|ds \leq C\int_{\R}N(s)e^{\pm \delta(t-s)}ds = Ce^{\pm \delta t}\int_{\R}N(s) e^{\mp \delta s}ds =: C'  e^{\pm \delta t}.
$$ 
Thus
we conclude that the Laplace transform $F(z)$ of $C^2$-smooth function $f(t), $$\ |f|_{C^2} < \infty,$ is well defined in the strip
$-\delta < \Re z < \delta$, where we have  
$$
|F(z)| \leq \frac{C_{pq}}{|z|^2},   \quad   p \leq \Re z \leq q,  \quad p,q \in (-\delta, \delta). 
$$
Hence, we can apply the Laplace transform to the equation 
$$
v(t)+ f(t) = \int_{\R}N(s)v(t-s)ds, 
$$
to obtain that
$$
V(z)= \frac{F(z)}{\mathcal{N}(z)-1}, \quad  -\delta < \Re z <0,
$$
where $\mathcal{N}(z): = \int_{\R}e^{-zs}N(s)ds$ of $N$ is a regular function 
in the strip $|\Re z| <\delta$.  Observe also that 
$$
\mathcal{N}(0)=1, \quad \mathcal{N}'(0) =0, \quad \mathcal{N}''(0) = \int_{\R}s^2N(s)ds >0. 
$$
Now, since $V(z)$ is analytical in the strip $\Pi := \{-\delta < \Re z <0\}$, the function  
$ F(z)/(\mathcal{N}(z)-1)$ has the same property in $\Pi$.  On the other hand, for an appropriate 
$\delta' \in (0,\delta)$ the quotient  $F(z)/(\mathcal{N}(z)-1)$ defines a meromorphic function 
in $\Pi' := \{-\delta < \Re z <\delta'\}$, with a unique singularity  (double pole) at $z=0$. Note that Corollary \ref{inn} as well as the last argument in the proof of Lemma \ref{tg}
are used at this stage. 
 Since the Laplace transform $V$ of $v \in C^2(\R)$ is integrable along each vertical line inside of $\Pi$, we may apply the inversion formula to get, 
for  arbitrarily fixed $c \in (-\delta,0), \  r \in (0,\delta')$, 
$$
v(t) = \frac{1}{2\pi i}\int_{c-i\cdot \infty}^{c+i\cdot \infty}\frac{e^{zt}F(z)}{\mathcal{N}(z)-1}dz=  {\rm Res}_{z=0}\frac{e^{zt}F(z)}{\mathcal{N}(z)-1}+
\frac{1}{2\pi i}\int_{r-i\cdot \infty}^{r+i\cdot \infty}\frac{e^{zt}F(z)}{\mathcal{N}(z)-1}dz. 
$$
Next, observe that if $f(t)\equiv 0$ then also $F(z)\equiv 0$ so that $v(t) \equiv 0$. Therefore the only case of the interest is when $f(s') >0$ at some $s' \in \R$ that implies $F(0) >0$. 
Now, in such a case, we have that 
$$
|\int_{r-i\cdot \infty}^{r+i\cdot \infty}\frac{e^{zt}F(z)}{\mathcal{N}(z)-1}dz|\leq  c_0e^{rt} \int_{\R}\frac{ds}{r^2+s^2} \leq c_1 e^{rt}, \  t \in \R, 
$$
while a direct calculation shows that 
$$ 
 {\rm Res}_{z=0}\frac{e^{zt}F(z)}{\mathcal{N}(z)-1} = \frac{2F(0)}{ \mathcal{N}''(0)}t + \frac{F'(0)}{ \mathcal{N}''(0)}- \frac{2F(0) \mathcal{N}'''(0)}{ 3(\mathcal{N}''(0))^2} =: At+B, \quad A>0. 
$$
In consequence, as $t \to -\infty$,  
$$
v(t) = At +B +O(e^{rt}),  \quad {\rm with} \ A, r >0, 
$$
which contradicts to the boundary condition $v(-\infty) =1$.   \hfill {}
\end{proof}
Now we are ready to prove our first uniqueness result: 
\begin{theorem} \label{mth1}Assume  {\bf {(SB)}} except $\gamma_\phi < \gamma_K $  as well as  {\bf {(EC$_{\gamma_\phi}$)}}    and suppose further that  $\chi(0)<0, \ \chi(\gamma_K-)\not= 0,$
\begin{equation} \label{gvt}
|g(u,\tau) - g(v,\tau)| \leq g'(0,\tau)|u-v|, \ u,v \geq 0.
\end{equation}
Then equation (\ref{17}) has at most one bounded positive  solution $\varphi, \  \varphi(-\infty) =0$. 
Furthermore, $\gamma_\phi$ coincides with the minimal  positive zero $\lambda_l$ of $\chi(z)$ and such a solution (if exists) has the following representation: 
$$
\varphi(t+m) = (a-t)^ke^{\lambda_l t} + e^{(\lambda_l + \delta)t}r(t),\  { \ with \,  continuous \ } r \in L^2(\R), 
$$
for some appropriate $a,m \in \R$, $\delta >0$. Here $k=0$ [respectively, $k=1$] if  $\lambda_l$ is a simple [respectively, double] root 
of $\chi(z)=0$. 
\end{theorem}
\begin{remark}\label{kf}
By Lemma \ref{nond},   $\chi(\gamma_K-)\not= 0$ yields $\gamma_\phi = \lambda_l< \gamma_K$.  
We  assume this stronger assumption instead of  $\gamma_\phi < \gamma_K$ since it is more easy to use. 
In the section of applications,  the condition  $\chi(\gamma_K-)\not= 0$ is slightly modified in order
to take into account the dependence of $\chi, \gamma_K$ on the wave velocity $c$.  Recall that we need $\gamma_\phi < \gamma_K$  to apply the bootstrap argument. 
\end{remark}
\begin{proof}  \underline{Step I: Asymptotic behavior at $-\infty$}. It is clear that equation  (\ref{17}) can be written as the linear inhomogeneous equation
\begin{equation} \label{UE}
\varphi(t)  = \int_Xd\mu\int_{\R}K(s,\tau)g'(0,\tau)\varphi(t-s)ds+
\mathcal{D}(t), \ t \in \R,
\end{equation}
where all integrals are converging and 
$$\mathcal{D}(t):= \int_Xd\mu\int_{\R}K(s,\tau)(g(\varphi(t-s),\tau)- g'(0,\tau)\varphi(t-s))ds \leq 0, \ t \in \R.
$$
Take $C(\tau), \sigma, \zeta(x)$ as in  {\bf {(SB)}}.    Observe that without restricting the generality, we can assume in {\bf {(SB)}} that $(1+\alpha)\gamma_\phi < \gamma_K. $
Since equation (\ref{17})  is translation invariant, we can suppose that $\varphi(t) < \sigma$ for $t \leq 0$. 
Applying the bilateral Laplace transform to (\ref{UE}), we obtain that 
$$
\chi(z)\Phi(z) = \mathbf{D}(z).
$$
We claim that, due to conditions {\bf {(SB)}} and {\bf{ (EC$_{\gamma_\phi}$)}}, function $\mathbf{D}$ is regular in the strip $\Pi_\alpha=\{z:\Re z  \in (0, (1+\alpha)\gamma_\phi)\}$.  Indeed, we have
$$
\mathbf{D}(x+iy) = \int_{R}e^{-iyt}[e^{-xt}\mathcal{D}(t)]dt. 
$$
Given $x: = \Re z  \in (0, (1+\alpha)\gamma_\phi)$, we choose  $x'$ sufficiently close from the left to $\gamma_\phi$ to satisfy
$
-x+(1+\alpha)x' >0.
$
Then
$$
|e^{-xt}\mathcal{D}(t)|\leq 
e^{-xt}\left[\int_XC(\tau)d\mu\int^{+\infty}_tK(s,\tau)C_{x'}^{1+\alpha}e^{(1+\alpha)x'(t-s)}ds+\right.$$
$$\left. +2|\varphi|_{\infty}\int_Xg'(0,\tau)d\mu\int_{-\infty}^tK(s,\tau)ds\right]\leq 
$$
$$e^{-xt}\left[e^{(1+\alpha)x't}C_{x'}^{1+\alpha}\zeta((1+\alpha)x')+2|\varphi|_{\infty}\int_Xg'(0,\tau)d\mu\int_{-\infty}^tK(s,\tau)ds\right]=:
$$
$$e^{-xt}\left[e^{(1+\alpha)x't}A_1+2|\varphi|_{\infty}\int_Xg'(0,\tau)d\mu\int_{-\infty}^tK(s,\tau)e^{-(1+\alpha)x's}e^{(1+\alpha)x's}ds\right]\leq 
$$
$$e^{(-x+(1+\alpha)x')t}\left[A_1+2|\varphi|_{\infty}(1-\chi((1+\alpha)x'))\right]
=:A_2e^{(-x+(1+\alpha)x')t}, \ t \in \R. 
$$

Since clearly $\mathcal{D}(t)$ is bounded on $\R$, the above calculation shows that 
$e^{-xt}\mathcal{D}(t)$ belongs to $L^k(\R)$, for each $k \in [1, \infty]$ once $x \in (0, (1+\alpha)\gamma_\phi)$.  As a consequence, for each such  $x$ the function 
$\mathbf{d}_x(y):=\mathbf{D}(x+i\cdot y )$ is bounded and square integrable on $\R$. 

By our assumptions, $\chi(z)$ is also regular in the domain $\Pi_\alpha$,
while 
$
\Phi(z) = {\mathbf{D}(z)}/\chi (z)
$
is regular in  $\Re z  \in (0, \gamma_\phi)$ and meromorphic in $\Pi_\alpha$. In virtue of Lemma \ref{tg}, we can suppose that 
$\Phi(z)$ has a unique singular point $\gamma_\phi$ in $\Pi_\alpha$ which is either simple or double pole. 

Now, for some $x'' \in (0,\gamma_\phi)$,  using the inversion theorem for the Fourier transform, we obtain that for an appropriate sequence of integers $N_j \to +\infty$ 
$$
\varphi(t) = \frac{1}{2\pi i}\lim_{j \to +\infty} \int_{x''-iN_j}^{x''+iN_j}\frac{e^{zt}\mathbf{D}(z)}{\chi (z)}dz
$$
almost everywhere on $\R$, e.g. see  \cite[p. 9-10]{FA}.  
Next, if $x \in (\gamma_\phi, (1+\alpha\gamma_\phi))$ then 
$$
\int_{x''-iN}^{x''+iN}\frac{e^{zt}\mathbf{D}(z)dz}{\chi (z)} = \left( 
 \int_{x-iN}^{x+iN}+\int_{x''-iN}^{x-iN}- \int_{x''+iN}^{x+iN}\right)\frac{e^{zt}\mathbf{D}(z)dz}{\chi (z)} - 2\pi i{\rm Res}_{z=\gamma_\phi}\frac{e^{zt}\mathbf{D}(z)}{\chi (z)}.
$$
Since, by Corollary \ref{inn}, 
$$\lim_{j \to +\infty}\max_{z \in [x''\pm iN_j, x\pm iN_j]}(|\mathbf{D}(z)|+|1-\chi (z)|)=0,$$
we  conclude  that, for each fixed $t \in \R$ 
$$
\lim_{j \to +\infty}\int_{x''\pm iN_j}^{x\pm iN_j}\frac{e^{zt}\mathbf{D}(z)}{\chi (z)}dz=0.
$$
Therefore
$$
\varphi(t) =  - {\rm Res}_{z=\gamma_\phi}\frac{e^{zt}\mathbf{D}(z)}{\chi (z)} + \frac{e^{xt}}{2\pi}\int_{\R}\frac{e^{iyt}\mathbf{d}_x(y)}{\chi(x+iy)}dy.
$$
It should be noted here that $\mathbf{D}(\gamma_\phi) < 0$ since otherwise $\mathcal{D}(t)\equiv 0$ implying $\chi(z)\Phi(z)=\mathbf{D}(z) \equiv 0$ so that $\Phi(z) \equiv 0$,  a contradiction. Since 
$$
{\rm Res}_{z=\gamma_\phi}\frac{e^{zt}\mathbf{D}(z)}{\chi (z)}  = \frac{e^{\gamma_\phi t}\mathbf{D}(\gamma_\phi)}{\chi' (\gamma_\phi)},  \quad {\rm if} \ \lambda_l < \lambda_r, 
$$
$$
{\rm Res}_{z=\gamma_\phi}\frac{e^{zt}\mathbf{D}(z)}{\chi (z)}  = \frac{2e^{\gamma_\phi t}}{\chi'' (\gamma_\phi)}  \left(t \mathbf{D}(\gamma_\phi) + \mathbf{D}'(\gamma_\phi)- \mathbf{D}(\gamma_\phi)\frac{\chi'''(\gamma_\phi)}{3\chi''(\gamma_\phi)}\right),  \quad {\rm if} \ \lambda_l = \lambda_r,  
$$
we get the desired representation. 

\underline{Step II: Uniqueness}.
By the contrary, suppose that $\varphi_1$ and $\varphi_2$ are different solutions of 
(\ref{17})  in the sence that $\varphi_1(t) \not\in \{ \varphi_2(t+s), \ s \in \R\}$. Due to Step I  we may suppose  that $\varphi_1, \varphi_2$
 have the same main parts of their asymptotic representations:
$$
\varphi_j(t) = (a_j-t)^ke^{\gamma_\phi t} + e^{(\gamma_\phi + \delta)t}r_j(t), \ r_j \in L^2(\R).  
$$
Therefore $\omega(t): = \varphi_2(t)-\varphi_1(t) = e^{(\gamma_\phi + \delta)t}r(t), \ t \in \R, \ r \in L^2(\R),$ in the case 
of $\lambda_l <\lambda_r$ and 
$\omega(t) =  (a_2-a_1)e^{\gamma_\phi t}+ e^{(\gamma_\phi + \delta)t}r(t), \ t \in \R, \ r \in L^2(\R),$  in the case 
of $\lambda_l =\lambda_r$.
Set 
$$
w(t):= \int^t_{t-1}|\omega(s)|ds,
$$ 
it is clear that $w\in C^1(\R)$ is bounded and has bounded derivative on $\R$, in fact, $0< |w'|_{\infty}, |w|_{\infty} \leq \max\{|\varphi_1|_\infty, |\varphi_2|_\infty\}$. Furthermore,  if $\lambda_l <\lambda_r$ then 
$$
w(t) = \int^t_{t-1}|e^{(\gamma_\phi + \delta)s}r(s)|ds \leq e^{(\gamma_\phi + \delta)t}\int^t_{t-1}|r(s)|ds\leq e^{(\gamma_\phi + \delta)t} \sqrt{\int^t_{t-1}r^2(s)ds}, 
$$ 
so that $w(t) = e^{(\gamma_\phi + \delta)t} o(1)$ at $t =-\infty$.  Now, if  $\lambda_l =\lambda_r$, we know that 
$$
\omega(t) =  ae^{\gamma_\phi t}+ e^{(\gamma_\phi + \delta)t}r(t), 
$$
where we can suppose that $a \geq 0$. 
Therefore 
$$
-e^{(\gamma_\phi + \delta)t}|r(t)| \leq |\omega (t)| - ae^{\gamma_\phi t} \leq e^{(\gamma_\phi + \delta)t}|r(t)|, 
$$
so that, in view of the above estimation of $w(t)$, we get
$$
|\omega(t)| =  ae^{\gamma_\phi t}+ e^{(\gamma_\phi + \delta)t}r_1(t),  \ {\rm with }\  |r_1(t)| \leq |r(t)|, 
$$
$$
w(t)= \int_{t-1}^t|\omega(s)|ds =  \frac{a(1- e^{-\gamma_\phi})}{\gamma_\phi}e^{\gamma_\phi t}+ e^{(\gamma_\phi + \delta)t}o(1), \ t \to - \infty.  
$$
We have the following: 
 $$ \omega(t)
=  \int_Xd\mu(\tau)\int_\R K(s,\tau)(g(\varphi_2(t-s),\tau)-g(\varphi_1(t-s),\tau))ds, 
$$
$$ |\omega(t)|
\leq \int_Xg'(0,\tau)d\mu(\tau)\int_\R K(s,\tau)|\omega(t-s)|ds,
$$
$$\int_{t-1}^t|\omega(u)|du
\leq  \int_Xg'(0,\tau)d\mu(\tau)\int_\R K(s,\tau)\int^t_{t-1}|\omega(u-s)|duds, 
$$
and, finally,  
\begin{equation} \label{pe}
w(t)
\leq  \int_Xg'(0,\tau)d\mu(\tau)\int_\R K(s,\tau)w(t-s)ds. 
\end{equation}
 \underline{Case I (noncritical)}.  If $\chi'(\lambda_l)\not= 0$, then $\chi(\gamma') >0$ for some $\gamma' \in (\gamma_\phi, \gamma_\phi+ \delta).$  After multiplying the both sides of (\ref{pe}) by $e^{-\gamma' t}$ and setting 
$
v(t):= w(t)e^{-\gamma' t}, 
$ 
we find that 
$$
v(t) \leq \int_\R \left( \int_Xg'(0,\tau)K(s,\tau)e^{-\gamma' s}d\mu(\tau)\right)v(t-s)ds. 
$$
Since $v(t)\geq 0$ and $v(\pm\infty)=0$, there exists a finite $t_m$ such that 
$$
v(t_m) = |v|_\infty = \max_{s \in \R} v(s). 
$$
But then 
$v(t_m) \leq \left( \int_Xg'(0,\tau)d\mu(\tau)\int_\R  K(s,\tau)e^{-\gamma' s}ds\right)v(t_m)$,
forcing $0=v(t_m)\equiv v(t)\equiv w(t)$ in view of $\chi(\gamma') >0$. 

 \underline{Case II (critical)}.  Now, if $\lambda_l=\lambda_r$,  we set  $
v(t):= w(t)e^{-\gamma_\phi t}, 
$  to  conclude analogously  that
$v(-\infty) = a(1-e^{-\gamma_\phi})/\gamma_\phi, \ v(+\infty) =0$, 
$$
v(t) \leq \int_\R \left( \int_Xg'(0,\tau)K(s,\tau)e^{-\gamma_\phi s}d\mu(\tau)\right)v(t-s)ds. 
$$
After normalizing if necessary, we can assume that 
$0\leq v(t) \leq 1= \sup_{s \in \R} v(s)$  for all $t\in \R$.  If $v(\hat t) =1$ for some finite  rightmost $\hat t$, then 
$$
1= v(\hat t) \leq \int_\R \left( \int_Xg'(0,\tau)K(s,\tau)e^{-\gamma_\phi s}d\mu(\tau)\right)v(\hat t-s)ds = :
$$
$$
\int_\R N(s)v(\hat t-s)ds
\leq  
 \int_Xg'(0,\tau)d\mu(\tau)\int_\R K(s,\tau)e^{-\gamma_\phi s}ds=1,
$$
which implies that $N(s)v(\hat t -s) = N(s)$ a.e. and $v(\hat t -s)=1$ for all $s$ such that $N(s)>0$.  Now, 
since $\int_{\R}N(s)ds=1, \ \int_{\R}sN(s)ds=0$, there is a subset of $\R_-$ of positive measure where 
$N(s) >0$. This means that $\hat t$ does not possesses  the property to be the rightmost point where $v(\hat t)=1$, a contradiction. 
Thus we have to analyze only the case when $a > 0$ and $0\leq v(t) < 1 = v(-\infty)$.  It is easy to check that in such a case, 
$v(t)$ and $N(s)$ meet all the conditions of Lemma \ref{trud}. In particular, since $\gamma_\phi < \gamma_K$, there exists $\delta >0$ such that
$$
\int_{\R}N(s)e^{xs}ds = 1-\chi(\gamma_\phi-x) < \infty \quad {\rm for \ all \ } |x| < \delta.  
$$
Hence, $v(t)\equiv 0$, a contradiction.  
\end{proof}
Let us consider now the situation  when the subtangential Lipschitz condition of Theorem \ref{mth1} is not satisfied.  In such a case, we prove the uniqueness under somewhat stronger 
hypotheses ({\bf{SB*}), ({\bf{EC*}}}): 
\begin{description}
\item[ {\bf{ (SB*)}}] Either one of the following conditions holds
$$
|g(u,\tau)-g(v,\tau)- g'(0,\tau)(u-v)| \leq C(\tau)|u-v|^{1+\alpha}, \ u,v \in (0,\sigma), 
$$
$$
|g'(u,\tau)- g'(0,\tau)| \leq C(\tau)u^{\alpha}, \ u \in (0,\sigma), 
$$
for some $\alpha, \sigma \in (0,1]$ and measurable $C(\tau) >0$ satisfying (\ref{ze}).  Furthermore, there exist  $\hat \epsilon \in (0,\gamma_\phi)$ and measurable 
$d_1(\tau)$ such that  
$$
0\leq K(s,\tau) \leq d_1(\tau)e^{\hat \epsilon s}, \ s \in \R. 
$$
\end{description}

\begin{description}
\item[ {\bf{ (EC*)}}]Either one of the following two assumptions is satisfied:
\begin{itemize}
\item[(i)]  Each solution of (\ref{17}) is $C^1$-smooth and if $\varphi_1, \varphi_2  \in C^1(\R)$ satisfy (\ref{17}) and the integral $\int_{\R}e^{-zs}(\varphi_2(s)-\varphi_1(s))ds$ converges  absolutely then the integral $\int_{\R}e^{-zs}(\varphi_2'(s)-\varphi_1'(s))ds$ also converges  absolutely. 
\item[(ii)]  There exists  $\delta_0 >0$ such that, for each $x \in (\lambda_{rK}-\delta_0, \lambda_{rK})$, it holds 
$$
0\leq K(s,\tau) \leq d_{2x}(\tau)e^{x s}, \ s \in \R, 
$$
for some $\mu-$measurable 
$d_{2x}(\tau)$.   
\end{itemize}
\end{description}
\begin{theorem} \label{mth2} Assume {\bf {(SB*)}}, {\bf{ (EC*)}} and suppose that
\begin{equation}\label{glc}
|g(u,\tau) - g(v,\tau)| \leq \lambda(\tau)|u-v|, \ u,v \geq 0, \tau \in X, 
\end{equation}
for some measurable $\lambda(\tau)$ different from $g'(0,\tau)$  and that the
function 
$$
\chi_1(z) =1 - \int_\R \int_XK(s,\tau)\lambda(\tau)d\mu(\tau)e^{-sz}ds 
$$
is well defined on $[0, \lambda_{rK})$. If, in addition,  $\lambda d_j \in L^1(X), \  j =1,2,$  $\chi(0) <0$ and
$\chi_1(m) \geq 0
$ for some $m \in (0, \lambda_{rK})$, 
then equation (\ref{17}) has at most one bounded positive  solution $\varphi, \  \varphi(-\infty) =0$. 
Furthermore, $\gamma_\phi$ coincides with the minimal  simple positive zero $\lambda_l$ of $\chi(z)$ and, for  appropriate $m \in \R$, $\delta >0$,  $$
\varphi(t+m) = e^{\lambda_l t} + e^{(\lambda_l + \delta)t}r(t),\  { \ with \,  continuous \ } r \in L^2(\R). 
$$
\end{theorem}
\begin{proof}  Using Lemma \ref{nond} and the above conditions,  we find that $\lambda_l = \gamma_\phi < m < \lambda_{rK}\leq \gamma_K$. Hence, due to Lemma \ref{eq},   the assumptions of the theorem  guarantee the  fulfillment of the hypotheses {\bf {(SB)}} and {\bf{ (EC$_{\gamma_\phi}$)}}.   Therefore all arguments of Step I in the proof of Theorem \ref{mth1} can be repeated  (with a unique change in the estimation of $e^{-xt}\mathcal{D}(t)$ where $g'(0,\tau), \chi$ is replaced with $\lambda(\tau), \chi_1$). 
Thus each pair $\varphi_1, \varphi_2$ of solutions of 
(\ref{17}) can be supposed to
 have the same main parts of their asymptotic representations:
$
\varphi_j(t) = e^{\lambda_l t} + e^{(\lambda_l + \delta)t}r_j(t), \ r_j \in L^2(\R).  
$ 
The further proof is  divided in  three steps. 

\underline{Step I}. Again, we  consider bounded function  $\omega(t): = \varphi_2(t)-\varphi_1(t) = e^{(\lambda_l + \delta)t}r(t),$ $ t \in \R, \ r \in L^2(\R)$. If $\Re z \in (0, \lambda_l+\delta)$, then   $\int_{\R}e^{-zs}\omega(s)ds$ converges absolutely  and from condition {\bf{ (EC*)}}(i)  we have
$$
|\omega(t)| = |\int_{-\infty}^t\omega'(s)ds| = |\int_{-\infty}^te^{xs}\omega'(s)e^{-xs}ds| \leq e^{xt}  
\int_{\R}e^{-xs}|\omega'(s)|ds=: C_x e^{xt}, 
$$
for all $x \in (0, \lambda_l+\delta)$ and $ \ t \in \R$.  Similarly,  we obtain from  {\bf {(SB*)}}, {\bf{ (EC*)}}(ii) that 
$$
|\omega(t)| = |\int_Xd\mu\int_{\R}K(s,\tau)\Big(g(\varphi_1(t-s),\tau)-g(\varphi_2(t-s),\tau)\Big)ds| \leq 
$$
$$
e^{xt}\int_X\lambda(\tau)d\mu\int_{\R}K(s,\tau)e^{-xs}e^{-x(t-s)}|\omega(t-s)|ds \leq 
$$
$$
e^{xt}\int_X\lambda(\tau)(d_1(\tau)+d_{2,\lambda_l+\delta}(\tau))d\mu\int_{\R}e^{-xs}|\omega(s)|ds, \ x \in (\hat\epsilon, \lambda_l+\delta), \ t \in \R.
$$
In each of these two cases,  for every   $x \in (\hat\epsilon, \lambda_l+\delta)$ there exists 
an  appropriate $C_x >0$ such that
$
|\omega(t)|\leq  C_x e^{xt},   \ t \in \R.
$
Set 
$$
\Gamma = \sup\{x \geq \lambda_l | \  \exists C_x : |\omega(t)|\leq  C_x e^{xt},   \ t \in \R \}, 
$$
we claim that $\Gamma \geq \lambda_{rK}$.  Indeed,  on the contrary, suppose that 
$\Gamma < \lambda_{rK}$ and let $x_0 \in (\hat\epsilon, \Gamma), \alpha >0, \gamma_0 \in (\hat\epsilon, \lambda_l)$ be such that $\{x_0(1+\alpha), \ 
x_0 + \alpha \gamma_0\}  \subset ( \Gamma, \lambda_{rK})$.  Let $x_*$ be the minimal  of these 
two numbers.  
We have that
\begin{equation} \label{UE2}
\omega(t)  = \int_Xd\mu\int_{\R}K(s,\tau)g'(0,\tau)\omega(t-s)ds+
\mathcal{E}(t), \ t \in \R,
\end{equation} with  bounded
$$
 \mathcal{E}(t):= \int_Xd\mu\int_{\R}K(s,\tau)\Big(g(\varphi_1(t-s),\tau)-g(\varphi_2(t-s),\tau)- g'(0,\tau)\omega(t-s)\Big)ds. $$
 Now, depending on assumptions chosen in  {\bf {(SB*)}}, we have either 
 $$
 |g(\varphi_1(s),\tau)-g(\varphi_2(s),\tau)- g'(0,\tau)\omega(s)| \leq C(\tau)|\omega(s)|^{1+\alpha} \leq 
 $$
 $$
 C(\tau)\min\{C_{x_0}^{1+\alpha}e^{x_0(1+\alpha)s},  ( |\varphi_1|_\infty+ |\varphi_2|_\infty)^{1+\alpha}\}
 \leq k_1C(\tau)e^{x_*s}, \ s \in \R, 
 $$
 or 
 $$
 |g(\varphi_1(s),\tau)-g(\varphi_2(s),\tau)- g'(0,\tau)\omega(s)| \leq C(\tau)|\omega(s)|(|\varphi_1(s)|+|\varphi_2(s)|)^\alpha \leq 
 $$
 $$
 k_2C(\tau)\min\{C_{x_0}e^{(x_0+\alpha\gamma_0)s},  (|\varphi_1|_\infty+ |\varphi_2|_\infty)^{1+\alpha}\}
 \leq k_3C(\tau)e^{x_*s}, \ s \in \R,
 $$
 where $k_i$ depend on $x_0,\gamma_0$ and $|\varphi_j|_\infty$ only. Hence, 
$$
|\mathcal{E}(t)|\leq  
 4e^{x_*t}\max\{|\varphi_1|_\infty,|\varphi_2|_\infty\}\int_X\lambda(\tau)d\mu\int_{-\infty}^t K(s,\tau)e^{-x_*s}ds
$$
$$
+ke^{x_*t}\int_XC(\tau)d\mu\int_{t}^{+\infty}K(s,\tau)e^{-x_*s}ds\leq $$
$$
 e^{x_*t}\Big (4\max\{|\varphi_1|_\infty,|\varphi_2|_\infty\} (1-\chi_1(x_*))+ k \zeta(x_*)\Big )=:Ae^{x_*t}, \ t \in \R. $$
 Therefore $e^{-xt}\mathcal{E}(t)$ belongs to $L^k(\R)$, for each $k \in [1, \infty]$ once $x \in (0, x_*)$. 
 Using Lemma \ref{tg}, we can repeat now the arguments of Step I  of Theorem \ref{mth1} (below the estimation of $|e^{-xt}\mathcal{D}(t)|$) to conclude that 
 $\omega(t)=e^{x t}r_x(t)\,\, \ t \in \R, \ r_x \in L^2(\R),$ for each $x  \in (\lambda_l, x_*)$. This implies the absolute 
 convergence of  $\int_{\R}e^{-xs}\omega(s)ds$  for every $x  \in (\lambda_l, x_*)$. But as we have seen at the 
 beginning of Step I, this  yields  $|\omega(s)| \leq B_xe^{xs}, \ s \in \R, \ x \in (\lambda_l, x^*)$ for  appropriate $B_x$. Therefore $\Gamma \geq x_* > \Gamma$, a contradiction.  In this way, we have proved that 
\begin{equation} \label{ome}
|\omega(s)| \leq B_xe^{xs}, \ s \in \R, \ x \in (\hat\varepsilon, \min\{\lambda_r, \gamma_K\}). 
\end{equation}
 
\underline{Step II}.  Suppose that $\chi_1(m) >0$ for some $m \in (0, \lambda_{rK})$,  it is clear that  $m > \lambda_l$ and 
$$
\kappa:=\int_\R \int_XK(s,\tau)\lambda(\tau)d\mu(\tau)e^{-sm}ds<1.$$
We now define $\bar\omega(t):=|\omega(t)|e^{-mt}\geq 0, \,t\in \R$. By (\ref{ome}),  we obtain that  $\bar\omega(\pm\infty)=0$ and $\bar\omega(t_m)= \max_{s \in \R}
\bar\omega(s)\geq 0$  for some $t_m\in \R$. Since
$$ \omega(t)
=  \int_Xd\mu(\tau)\int_\R K(s,\tau)(g(\varphi_2(t-s),\tau)-g(\varphi_1(t-s),\tau))ds, 
$$
we have
$$
\bar\omega(t_m)=|\omega(t_m)|e^{-mt_m}
\leq \int_X\lambda(\tau)d\mu(\tau)\int_\R K(s,\tau)e^{-ms}|\omega(t_m-s)|e^{-m(t_m-s)}ds
\leq
$$
$$
 \bar\omega(t_m)\int_X\lambda(\tau)d\mu(\tau)\int_\R K(s,\tau)e^{-ms}ds=\bar\omega(t_m)\kappa.
$$
Hence, 
 $\bar\omega(\tau)=0$ and the uniqueness  follows.
 
 \underline{Step III}. Suppose now that $\chi_1(m)=\max_{s \in (0, \lambda_{rK})}\chi_1(s)= 0$.  Then additionally  $\chi_1'(m)=0$. Since $\lambda(\tau)$ is different from $g'(0,\tau)$, we have also that $\lambda_l <m$. Furthermore, $\bar\omega(t):=|\omega(t)|e^{-mt}\geq 0, \,t\in \R$ has the same properties as in Step II: 
$ \bar\omega(\pm\infty)=0$, \ 
 $\bar\omega(t_m)= \max_{s \in \R}
\bar\omega(s)\geq0$  for some $t_m\in \R$ and 
$$
\bar\omega(t)\leq \int_\R \left( \int_X K(s,\tau)\lambda(\tau)e^{-m s}d\mu(\tau)\right)\bar\omega(t-s)ds. 
$$
After normalizing, we may assume that 
$0\leq \bar\omega(t) \leq 1=\bar\omega(t_m) =1$, $ t\in \R$, for some finite  rightmost $ t_m$. Then 
$$
1\leq \int_\R N_\lambda(s)\bar\omega(t_m-s)ds\leq \int_\R N_\lambda(s)ds=1,
$$
where $N_\lambda(s):=\int_X K(s,\tau)\lambda(\tau)e^{-m s}d\mu(\tau)$. This implies  that $N_\lambda(s)\bar\omega(t_m -s) $ $= N_\lambda(s)$ a.e. and $\bar\omega(\hat t -s)=1$ for all $s$ such that $N_\lambda(s)>0$.  Now, 
since $\int_{\R}N_\lambda(s)ds=1, \ \int_{\R}sN_\lambda(s)ds=0$, there is a subset of $\R_-$ of positive measure where 
$N_\lambda(s) >0$. This means that $t_m$ does not possesses  the property to be the rightmost point where $\bar\omega(t_m)=1$, a contradiction. In consequence,   $\bar\omega(t)\equiv 0 $ that proves the uniqueness. 
\end{proof}
\begin{remark} It is enlightening to compare Theorem \ref{mth2} and Theorem 2 in \cite{KS} where somewhat similar ideas were exploited.  Indeed, from pure analytical estimations,  without the use of asymptotic representations of solutions and without using the properties of $\chi$ indicated in Lemma \ref{tg}, Schumacher deduced that $\Gamma \geq \lambda_{rK}$ (under assumptions made in \cite{KS}).  In any case,  monotonicity 
restrictions on the convolution term in \cite{KS} do not allow consider  various interesting models (cf. Sections 6.3-6.4 below). 
\end{remark}

\section{Applications}

In this section, Theorems \ref{mth1}  and \ref{mth2}  are applied to several models which can be written  as  (\ref{17}).  This allows to  improve or complement the uniqueness results in \cite{atv,CC,CDM,dk,fwz,tz}.
Everywhere in  this section we assume that locally Lipschitzian  $g: \R_+ \to \R_+, \ g(0)=0,$ is differentiable at $0$ with  $g'(0) >0$.
\vspace{-3mm}
\subsection{A nonlocal integro-differential equation \cite{CC,CDM,CDu,Ga,KPP,KS}} \label{ssc}
Consider  the equation
\begin{equation}\label{cc}
u_t=J*u -u+g(u),
\end{equation}
where $J\geq 0$, $\int_\R Jds >0$.  Let  $\gamma^\#$ denote an extended positive real number such that $\int_\R J(s)e^{-zs}ds$ is convergent  for $z\in[0,\gamma^\#)$ and is divergent when $z > \gamma^\#$. As it can be easily deduced from Theorem \ref{thm:expnv}, the existence of such $\gamma^\#$ is automatically assured by 
the existence of positive semi-wavefronts $u(t,x)=\phi(x+ct), \ \phi(-\infty)=0 $ to (\ref{cc}).  
Traveling wave profile $\phi$ solves
\begin{equation}\label{cc2}
c\phi'=J*\phi -\phi+g(\phi).
\end{equation}
In order to replace condition 
(\ref{lcs}) with more weak   
\begin{equation}\label{lae}
g'(s) \leq g'(0)  \ {\rm a.e. \  on} \  \R_+, 
\end{equation}
we use the following trick.  Set $g_\beta(s) = g(s) + \beta s$ for some positive $\beta$.  
We claim that  $\beta$ can be chosen in such a way that $g_\beta$ satisfies the Lipshitz condition with a constant $g_\beta'(0)=\beta +  g'(0)$.  First observe that our proof of uniqueness compares two different solutions $\phi_1, \phi_2$. Since they are uniformly bounded by some positive $M>0$, we can restrict our attention to a finite interval $[0,M]$ where $g$ is  globally Lipschitzian.  But then there exists $\beta>0$  such that $g'(0)\geq g'(s)\geq -2\beta - g'(0)$ almost everywhere on  $[0,M]$. In consequence, we get the necessary estimation 
$$- g'(0) - \beta \leq g'_\beta(s)= \beta + g'(s) \leq \beta + g'(0) \ {\rm a.e. \  on} \  [0,M].$$
Hence, instead of (\ref{cc2}) we will consider  
\begin{equation}\label{cc2b}
c\phi'=J*\phi -(1+\beta)\phi+g_\beta(\phi). 
\end{equation}
Let us suppose that $c>0$ (the case $c<0$ is similar). Since $\phi$ is non-negative and bounded, it should satisfy 
$$
\phi(t)=\frac{1}{c}\int_{-\infty}^{t} e^{-(t-s)(1+\beta)/c}\big(J*\phi (s)+g_\beta(\phi(s))\Big )ds=
$$
\begin{equation}\label{cc3}
k*(J*\phi) (t)+k*g_\beta(\phi) (t)
=(k*J)*\phi (t)+k*g_\beta(\phi) (t),
\end{equation}
where $k(s)=c^{-1}e^{-s(1+\beta)/c},\,s\geq 0$ and $k=0$ if $s<0$. Thus,  equation (\ref{cc3}) can be written as  (\ref{tweabnv}),  with $X=\{\tau_1,\tau_2\}$ and 
$$K(s,\tau)=\left\{\begin{array}{cc} k*J (s), & \tau=\tau_1, \\ k(s), & \tau=\tau_2, \end{array}\right.\quad  g(s,\tau)=\left\{\begin{array}{cc}  s, & \tau=\tau_1, \\ g_\beta(s),& \tau=\tau_2.\end{array}\right.$$ 
Finally, independently on the sign of $c$, we find that 
 $$
\chi(z,c)=1-\int_\R K(s,\tau_1)e^{-zs}ds-(g'(0)+\beta)\int_\R K(s,\tau_2)e^{-zs}ds=
$$
$$
 1-\frac{1}{1+\beta+ cz}\int_\R J(s)e^{-zs}ds-\frac{g'(0)+\beta}{1+\beta+cz} =: \frac{\tilde \chi(z,c)}{1+\beta+cz}.
$$
Let $c_*$ be the minimal value of $c$ for which  
$$\tilde\chi(z,c):= 1-g'(0) +cz - \int_\R J(s)e^{-sz}ds$$
has at least one positive zero.  It is easy to see that 
$$
c_* = \inf_{z>0}\frac 1z\left\{-1+g'(0) + \int_\R J(s)e^{-sz}ds\right\}
$$
can be positive, negative (in these cases $\inf$ can be replaced with $\min$) or zero. 
 By Theorem \ref{gg},  $c \geq c_*$ for each admissible wave speed $c$. 
The next result is a direct consequence of Theorem \ref{mth1}.
\begin{theorem} \label{th5}
Suppose (\ref{lae})  together with $1-\int_\R J(s)ds<g'(0)$ and
 \begin{equation}\label{**}
 |g(u)- g'(0)u| \leq Cu^{1+\alpha}, \ u,v \in (0,\sigma) \,\, {\rm \ for\ some \ }\,\, \alpha,\sigma\in (0,1],
 \end{equation}
Then equation (\ref{cc2}) has at most one bounded positive  solution $\varphi, \ \varphi(-\infty) =0$, for each $c\not=0$ (if $\tilde \chi(\gamma^\#-,c_*)\not=0$) or for each $c \not= 0, c_*$ (if $\tilde \chi(\gamma^\#-,c_*)=0$). 
\end{theorem}
\begin{proof} Suppose that $c>0$ (the case $c<0$ is similar). 
We only have to check the assumptions {\bf {(EC$_{\gamma_\phi}$)}}, {\bf {(SB)}} except $\gamma_\phi(c) < \gamma_K(c)$,  $\chi(0,c)<0$ and  $\chi(\gamma_K-,c)\not=0$
 of Theorem \ref{mth1}. \\
 \underline{Step I}. It is clear that $g(\cdot,\tau)$ satisfies (\ref{gvt}), where $g'(0,\tau_1)=1,\, g'(0,\tau_2)=g'(0)+\beta$.  Moreover, we have 
$
|g(u,\tau)- g'(0,\tau)u| \leq C(\tau)u^{1+\alpha}, \ u,v \in (0,\sigma)
$, where $C(\tau)=0$ if $\tau=\tau_1$ and $C(\tau)=C$ if $\tau=\tau_2$.\\
 \underline{Step II}. For each
  $z>-\frac{1+\beta}{c}$  we have $\int_\R k(s)e^{-zs}ds=\frac{1}{1+\beta +cz}<+\infty$ so that 
$\gamma_K(c)=\gamma^{\#}$ because of  
 $
\int_\R k*J(s)e^{-zs}ds=\int_\R J(s)e^{-zs}ds/(1+\beta+cz).
 $ (Observe here that $\gamma_K(c)=\min\{\gamma^{\#}, -(1+\beta)/c\}$ if $c <0$. However, if $\gamma_K(c)=-(1+\beta)/c$ then $\chi(\gamma_K(c),c)= \infty$ so that $\gamma_\phi(c) < \gamma_K(c)$ due to Corollary \ref{cfi}).\\ 
 \underline{Step III}. If  $\varphi$  solves (\ref{cc2}), then $\varphi\in C^1(\R)$ and for each $0<z<\gamma _\phi$ we obtain
 $$
 c\int_{\R}e^{-zs}|\varphi'(s)|ds\leq \int_{\R}e^{-zs}J*\varphi(s)ds+\int_{\R}e^{-zs}\varphi(s)ds +\int_{\R}e^{-zs}g(\varphi(s))ds
\leq 
$$
$$
(\int_{\R}e^{-zs}J(s)ds +1 +g'(0))\int_{\R}e^{-zs}\varphi(s)ds<+\infty.
$$
Thus, by Lemma  \ref{eq}, condition {\bf {(EC$_{\gamma_\phi}$)}} is satisfied. 
\\
 \underline{Step IV}. We have  $\chi(0,c)=(1-\int_\R J(s)ds-g'(0))/(1+\beta)<0$. Now, if $\gamma^\#<+\infty$, then $\tilde \chi(\gamma^\#-,c_*)\not=0$ implies that $\chi(\gamma^\#-,c_*)\not=0$ and $\gamma_\phi(c_*)=\lambda_l(c_*) < \gamma^\#$.  Since $\chi(z,c)$ is strictly increasing in $c$ for each fixed $z >0$, function $\lambda_l(c)$ is strictly decreasing. Hence $\gamma_\phi(c)=\lambda_l(c) < \gamma^\# $ for each $c \geq c_*$. Similar considerations shows that  $\gamma_\phi(c) < \gamma^\# $ for each $c > c_*$ if $\chi(\gamma^\#-,c_*)=0$. Finally, in the case $\gamma^\#=+\infty$ we have that $\chi(+\infty,c) \in \{1, -\infty\} \not\ni 0$,  so that  $\chi(\gamma_K-,c)\not=0$ holds automatically. \end{proof}

\begin{remark} \label{ruma}
Our approach allows to remove several restrictions on $J$  and $g$ assumed in the Carr and Chmaj uniqueness  result \cite[Theorem 2.1]{CC}.  In the cited work   $g$ is supposed to satisfy (\ref{lcs}) and $J$  to be an even compactly supported function with $\int_\R Jds=1$. These properties were  essential in the proof of Theorem 2.1 in \cite{CC} even though  (\ref{lcs}) was not mentioned explicitly there.  Similarly,  conditions $J\in C^1(\R)$, $J(a)>0, J(b)>0$ for some $a<0<b$, and of $J$  compactly supported were  used by Coville {\it et al.} It was   assumed  in \cite{CDM} that $g'(0)g'(1) <0$ together with 
 $g(u)/u \leq g'(0), \ u > 0,$ instead of more restrictive $g'(u)\leq g'(0), \ u >0$.  See also \cite{CDM} for non-uniqueness of stationary traveling fronts ($c=0$).  Next,  Schumacher \cite{KS}, using 
a comparison method for differential inequalities combined with a Nagumo-point argument,   
 established uniqueness of {\it regular and non-critical}  semi-wavefronts to equation (\ref{cc}) for general $J$ and $g$ satisfying  (\ref{lae}).   The trick allowing to weaken the Lipschitz restriction (\ref{lcs}) is due to Thieme and Zhao \cite{tz}   (as far as we know).  Usually it was applied under reversed inequality $f'(s) \geq f'(0)$ to the second (damping) term of equation, e.g. see also \cite{fz} and Section 6.3 for further generalizations.  Here we show that this trick shows to be useful also in the case of birth functions.  We would like to note that Theorem \ref{th5} remains true if we introduce a small delay $h >0$ in the term $g(\varphi(t-h))$. Indeed, in such a case it suffices to replace $k(s)$ with a positive fundamental solution $v(s)$
 of the scalar delayed equation $cv'(s) = - v(s) - \beta v(s-h)$.

 \end{remark}

\subsection{Nonlocal lattice equations \cite{CG,fwz,GW,MaZ,zhh}} 
Now we consider semi-wavefronts $w_j(t)=u(j+ct), u(-\infty)=0$, of the nonlocal lattice equation
$$
\frac{dw_j(t)}{dt}= D[w_{j+k}(t)-w_j(t)]-dw_j(t)+ \sum_{k\in \Z}\beta(j-k)g(w_k(t-r)),  \ j \in \Z, $$
where $\beta(k)\geq 0, \ \sum_{k\in \Z}\beta(k)=1$. 
Let  $\gamma^\#$ be an extended positive real number such that $\sum_{k\in \Z}\beta(k)e^{-zk}$ converges when $z\in[0,\gamma^\#)$ and is divergent when $z > \gamma^\#$.  By Cauchy-Adamard formula, 
$\gamma^\# = - \limsup_{k \to +\infty} k^{-1}\ln\beta(-k)$, where by convention $\ln(0)=-\infty$.
The wave profile $u$ satisfies 
\begin{equation}\label{sy1}
cu'(x)= D[u(x+1)+u(x-1) -2u(x)] - du(x) + \sum_{k\in \Z}\beta(k)g(u(x-k-cr)).
\end{equation}
Again we take $c>0$ for simplicity. Since  $u$ is bounded, we find that 
$$
u(t) = \frac1c \int^t_{-\infty}e^{-\frac{2D+d}{c}(t-s)}\left[Du(s+1)+Du(s-1)  + \sum_{k\in \Z}\beta(k)g(u(s-k-cr))\right]ds$$
\begin{equation}\label{sy2}\nonumber
=D(H_{-1}+H_1)*u(t)+ \sum_{k\in \Z}\beta(k)H_{k+cr}*g(u)(t),
\end{equation}
where 
$$H_\tau(t)=\left\{\begin{array}{cc} \frac 1c e^{-\frac{2D+d}{c}(t-\tau)},&t\geq \tau,  \\ 0,& t<\tau. \end{array}\right. $$
Thus (\ref{sy2}) can be written as (\ref{17}),  with $X=\{\tau_1,\tau_2\}$  and 
$$K(s,\tau)=\left\{\begin{array}{cc} D(H_{-1}(s)+H_1(s)), & \tau=\tau_1, \\  \sum_{k\in \Z}\beta(k)H_{k+cr}(s),& \tau=\tau_2, \end{array}\right. \quad  g(s,\tau)=\left\{\begin{array}{cc}  s, & \tau=\tau_1, \\ g(s),& \tau=\tau_2.\end{array}\right.$$ 
Next,  $\chi(z,c)=1- \int_\R K(s,\tau_1)e^{-sz}ds-g'(0)\int_\R K(s,\tau_2)e^{-sz}ds=$
$$
1-\frac{2D\cosh(z)}{2D+d+cz}-\frac{g'(0)e^{-crz}}{2D+d+cz}\sum_{k\in \Z}\beta(k)e^{-kz} =:\frac{\tilde\chi(z,c)}{2D+d+cz}.
$$
Let $c_*$ be the minimal value of $c$ for which  
$$\tilde\chi(z,c):= d+2D+cz -D(e^z +e^{-z})-g'(0)e^{-crz}\sum_{k\in \Z}\beta(k)e^{-kz}$$
has at least one positive zero.  It is easily  seen that $c_*$ is well defined and is finite. 
 By Theorem \ref{gg},  $c \geq c_*$ for each admissible wave speed $c$. 
 
 We are ready to  apply our uniqueness results to (\ref{sy1}).
\begin{theorem} \label{lat}
Suppose that $g$ satisfies (\ref{lcs}), (\ref{**})  and $g'(0)>d$.
Then equation (\ref{sy1}) has at most one bounded positive  solution $u, \ u(-\infty) =0$, for each $c\not=0$ (if $\tilde \chi(\gamma^\#-,c_*)\not=0$) or for each $c \not= 0, c_*$ (if $\tilde \chi(\gamma^\#-,c_*)=0$). 
\end{theorem}
\begin{proof}
\underline{Step I}. Obviously,  $g(\cdot,\tau)$ verifies (\ref{lcs}) with   $g'(0,\tau_1)=1$ and  $g'(0,\tau_2)=g'(0)$.  Moreover, we have 
$
|g(u,\tau)- g'(0,\tau)u| \leq C(\tau)u^{1+\alpha}, \ u,v \in (0,\sigma)
$, where $C(\tau_1)=0$  and $C(\tau_2)=C$.\\
 \underline{Step II}. If $0<z<\gamma_\#$, we get 
 $$
 \int_{\R\times X}K(s,\tau)e^{-zs}dsd\mu=\int_\R D(H_{-1}(s)+H_1(s))e^{-zs}ds+$$
 $$
 \int_\R \sum_{k\in \Z}\beta(k)H_{k+cr}(s)e^{-zs}ds =
\frac{2D\cosh(z)}{2D+d+cz}+\frac{e^{-cr}}{2D+d+cz}\sum_{k\in \Z}\beta(k)e^{-kz}.
$$
Therefore $\gamma_K=\gamma_\#$ (if $c >0$) and $\gamma_K=\min\{\gamma_\#, -(2D+d)/c\}$ (if $c <0$). \\
 \underline{Step III}. If  $u$  solves (\ref{sy1}) with $c>0$, then  for each $0<z<\gamma _\phi$ we obtain
 $$
c\int_\R |u'(s)|e^{-zs}ds \leq D\int_\R (u(s+1)+u(s-1) +2u(s)) e^{-zs}ds+
$$
$$d\int_\R u(s)e^{-zs}ds + g'(0)\sum_{k\in \Z}\beta(k)\int_\R u(s-k-cr)e^{-zs}ds=$$
$$
\left(2D(\cosh(z)+1) +d +
 g'(0)e^{-zcr}\sum_{k\in \Z}\beta(k)e^{-zk}\right)\int_\R u(s)e^{-zs}ds<+\infty.
$$
Thus, by Lemma  \ref{eq}, condition {\bf {(EC$_{\gamma_\phi}$)}} is satisfied.\\
 \underline{Step IV }. We have $\chi(0)=(d-g'(0))/(2D+d)<0$.  The proof of 
 $\gamma_\phi(c) < \gamma^\#$ is the same as in Step IV of the previous section and is omitted. 
\end{proof}
\begin{remark}
Our approach allows to improve the uniqueness results of \cite[Theorem 3.1]{fwz}, where additional conditions    $\beta(k)=\beta(-k)$ and $\chi(\gamma_K-)= -\infty$ are assumed.  Moreover, \cite[Theorem 3.1]{fwz} does not establish the uniqueness of the minimal wave. Similarly to  Section \ref{ssc}, condition (\ref{lcs}) in Theorem \ref{lat} can be replaced with more weak (\ref{lae}) if the nonlinear term  is local and non-delayed. See  \cite{GW}, where a  local and non-delayed variant of (\ref{sy1}) was considered. Similarly to \cite{Co,CDM} and under the same conditions on $g$ as in \cite{CDM},  Guo and Wu prove their uniqueness result  \cite[Theorem 2]{GW} by means of the comparison argument. To establish the uniqueness 
in  the degenerate case $(g'(0)-d)(g'(1)-d) = 0$ (cf. Remark \ref{ruma}), about which is the main concern of \cite{Chen},   Chen {\it et al.}  developed
new interesting tools (magnification, compression, blow-up techniques, modified sliding method). Finally, we  mention Ma and Zou uniqueness result from \cite{MaZ}, where a local version of   (\ref{sy1})
is investigated.  The Lipschitz condition (\ref{lcs}) is not required in \cite{MaZ}, it is supposed instead that $g'(s) \geq 0, \ g(s)/s \leq g'(0), \ s >0$. 

\end{remark}

\subsection{Nonlocal reaction-diffusion equation \cite{fz,gouss,MLLS,tz,TAT}}\label{nrd}
Here, we consider positive semi-wavefronts  $u(t,x)$ $= \phi(x +ct), \  \phi(-\infty) =0$, for  non-local delayed reaction-diffusion
equations  
\begin{equation}\label{r-d}
u_t(t,x) =  u_{xx}(t,x)- f(u(t,x)) +
 \int_{\R}k(w)g(u(t-h,x-w))dw,\,h>0,
 \end{equation}
 where $f\in C^1(\R_+,\R_+), f(0)=0,$ is strictly increasing and $k\geq 0, \  \int_{\R}kds=1,$ can be  asymmetric (see \cite{TAT} for  further details
concerning wave solutions in the presence of asymmetric non-local interaction). Let  $\gamma^\#>0$ denote an extended positive real number such that $\int_\R k(s)e^{-zs}ds$ converges when $z\in[0,\gamma^\#)$ and diverges if $z > \gamma^\#$. 
 It is clear that profile $\phi$  must satisfy
\begin{equation}\label{i2}
 y''(t) - cy'(t)-f(y(t))+ \int_{\R}k(s)g\left(y(t-ch-s)\right)ds=0, \ t \in \R.
\end{equation}
Equation (\ref{i2}) can be written as 
$$
 y''(t) - cy'(t)-\beta y(t)+f_\beta(y(t))+  \int_{\R}k_h(w)g(y(t-
w))dw=0,\, t \in \R,
$$
where  $k_h(w) = k(w-ch)$ and $f_\beta(s)=\beta s- f(s)$ for some  $\beta>0$.  

Again, without restricting the generality, we may suppose that  $f_\beta$ is a Lipshitzian function  with Lip$f_\beta=\beta - \inf_{s \geq 0} f'(s)$.  Indeed,  our proof of uniqueness compares two solutions $\phi_1, \phi_2$. Since they are uniformly bounded by some positive $M>0$, we can restrict our attention to a finite interval $[0,M]$. Let  $\beta>f'(0)$  be such that $f_\beta(s)=\beta s- f(s)\geq 0$ for all $s\in [0,M]$ and 
$\max _{s\in[0,M]}f'(s)\leq 2\beta-\inf_{s\geq 0}f'(s).$
But then 
$$
\left|\frac{f_\beta(s_2)-f_\beta(s_1)}{s_2-s_1}\right|\leq  \Big(\beta-\inf_{s\geq 0}f'(s)\Big),\quad s_1,s_2\in [0,M].
$$
Next, it is easy to see that the wave profile $\phi$ solves the equation
$$
\phi(t)=\frac{1}{\sigma(c)}\left(
\int_{-\infty}^te^{\nu(t-s)}(\mathcal G \phi)(s)ds
+\int_t^{+\infty}e^{\mu(t-s)}(\mathcal G\phi)(s)ds\right), 
$$
where  $\sigma(c)=\sqrt{c^2+4\beta}$,  $\nu<0<\mu$ are the roots of  $z^2-cz-\beta=0$
 and  $(\mathcal G\phi)(t):=\displaystyle \int_{\R}k_h(s)g(\phi(t-
s))ds+f_\beta(\phi(t))$. Equivalently, 
$$ 
\phi(t)=(\mathcal K*k_h)*g(\phi)(t)+\mathcal K* f_\beta(\phi)(t),$$
where
$$ \mathcal K(s)=\sigma^{-1}(c)\left\{\begin{array}{cc} e^{\nu s}, & s\geq 0, \\ e^{\mu s }, & s<0. \end{array}\right.
$$
We can invoke now Theorems \ref{mth1},  \ref{mth2}  where $X=\{\tau_1,\tau_2\}$ and 
$$K(s,\tau)=\left\{\begin{array}{cc} (\mathcal K*k_h)(s),& \tau=\tau_1, \\  \mathcal K,& \tau=\tau_2, \end{array}\right.\quad  g(s,\tau)=\left\{\begin{array}{cc}  g(s), & \tau=\tau_1, \\ f_\beta(s),& \tau=\tau_2.\end{array}\right.$$ 
Observe that  $g(\cdot,\tau)$  meets (\ref{glc}) with   $\lambda(\tau_1)=g'(0)$, $\lambda(\tau_2)=\beta-\inf_{s\geq 0}f'(s)$. If   $f'(0)\leq f'(v)$ for all $v\geq 0$, as in \cite{tz}, then $\beta-\inf_{s\geq 0}f'(s)=\beta-f'(0) = f'_\beta(0)$. We have also that 
$$
\chi_1(z,c) = 1-g'(0)\int_\R K(s,\tau_1)e^{-sx}ds-(\beta-\inf_{s\geq 0}f'(s))\int_\R K(s,\tau_2)e^{-sx}ds=$$
$$
 1-\frac{\beta-\inf_{s\geq 0}f'(s)}{\beta+cz-z^2}-\frac{g'(0)e^{-zch}}{\beta+cz-z^2}\int_\R k(s)e^{-zs}ds=: \frac{\tilde\chi_1(z)}{\beta+cz-z^2}.
$$
We  see that  $\gamma_K=\min\{\mu, \gamma^\#\}$ so that  $\gamma_\phi < \mu$. 
Let $c_\star$ be the minimal value of $c$ for which  
$$\tilde\chi_1(z,c):= cz - z^2 +\inf_{s\geq 0}f'(s)-g'(0)e^{-zch}\int_\R k(s)e^{-zs}ds$$
 has at least one positive zero.  This value is finite, well defined and does not depend on $\beta$.  We will write 
 $c_*$ instead of $c_\star$  in the special case  when $f'(0)\leq f'(v)$ for all $v\geq 0$. In such a case, we have  $f'(0)= \inf_{s\geq 0}f'(s)$ and therefore $\chi_1 = \chi$.  By Theorem \ref{gg},  $c \geq c_*$ for each admissible wave speed $c$. 
 \begin{theorem}\label{tz1}
Suppose $g$ satisfies (\ref{lcs}), $f\in C^1(\R_+, \R_+)$ is strictly increasing, and $g,f \in C^{1,\alpha}$ in some neighborhood of $0$, and $g(0)=f(0) =0,\ g'(0)>f'(0)$.  Then equation (\ref{r-d}) has at most one positive  semi-wavefront  $u(t,x)$ $= \phi(x +ct), \  \phi(-\infty) =0,$ for each $c\geq c_\star$ (if $\tilde \chi(\gamma^\#-,c_\star)\not=0$) or for each $c >c_\star$ (if $\tilde \chi(\gamma^\#-,c_\star)=0$). 
\end{theorem}
\begin{proof}
Observe that  $\beta \chi(0)=f'(0)-g'(0)<0,$  and 
$
\chi_1(\gamma^\#-,c_\star) \not= 0$ if $\tilde\chi_1(\gamma^\#-,c_\star) \not=0$. 
{\it First let} $c\geq c_\star > c_*$, then  $\chi_1(x,c) < \chi (x,c)$ so that   $\chi_1(m,c) = 0$ for some  $m \in (0, \lambda_{rK}]$.  It is clear that  $m = \lambda_{rK}$ if and only if $m= \gamma^\#$.  Since $\chi_1(z,c)$ is strictly increasing in $c$ for each fixed positive $z$, this implies that $c=c_\star$ and $\chi_1(\gamma^\#-,c_\star)=0$.  Consequently, $m \in (0, \lambda_{rK})$  for each $c\geq c_\star$ (if $\tilde \chi(\gamma^\#-,c_\star)\not=0$) or for each $c >c_\star$ (if $\tilde \chi(\gamma^\#-,c_\star)=0$).  \\
{\it Next, if}  $c_\star = c_*$ then $\chi_1=\chi$ and 
 the inequality $\chi(\gamma^\#-,c_*)\not=0$ guarantees that $\lambda_l(c_*) =\gamma_\phi(c_*) < \gamma^\#$ for $c=c_*$.  If $c>c_*$ then we have  again $\lambda_l(c)= \gamma_\phi(c)  <  \lambda_l(c_*)  < \gamma^\#$ because $\lambda_l(c)$ is monotone decreasing in $c$. \\
\underline{Step I}.  
 Since $|f_\beta'(0)u-f_\beta(u)|=|f'(0)u-f(u)|$, for an appropriate $C, \sigma$, it holds
$|g(u,\tau)- g'(0,\tau)u| \leq C(\tau)u^{1+\alpha}, \ u \in (0,\sigma)$. \\
\underline{Step II}. We claim that for each $x \in (0, \gamma_K)$ and some $d_j(x)$ it holds 
$$
0\leq K(s,\tau_j)\leq d_j(x)e^{x s}, s \in \R. 
$$
Indeed, if $j=2$,  we can even take $x= \mu, \ d_2 = 1/\sigma(c)$. Next,  we have
$$
K(t,\tau_1)=   \frac{1}{\sigma(c)}\left[\int_{t-ch}^{+\infty}e^{\mu (t-ch-v)}k(v)dv+ \int^{t-ch}_{-\infty}e^{\nu (t-ch-v)}k(v)dv \right ] \leq 
$$
$$
 \leq \frac{e^{-xch}}{\sigma(c)}\left[\int_{\R}e^{-xv}k(v)dv \right] e^{xt}.   
$$
Since $ \lambda_{rK}\leq \gamma_K = \min\{\gamma^\#,\mu\}$, the exponential estimations of $K$ in  {\bf {(SB*)}}, {\bf{ (EC*)}}(ii) are verified.  This observation completes the 
proof of  the theorem.  
\end{proof}
\begin{remark}
Theorem \ref{tz1} improves  \cite[Theorem 4.3]{tz}, where the uniqueness was established under assumption that either  $f(s)= f'(0)s$ or $g(s)=g'(0)s$ and $K$ is the Gaussian kernel.  Moreover, \cite[Theorem 4.3]{tz} does not consider the minimal waves.  See also \cite{Ma,TAT} and references therein about the existence of semi-wavefronts in (\ref{r-d})  and  its limit form (\ref{sv}) studied below. 
\end{remark}
\subsection{Uniqueness of fast traveling fronts in delayed  equations}
Here we study positive semi-wavefronts $u(t,x) = \phi(x
+ct), \phi(-\infty)=0,$ to
\begin{equation}\label{sv}
u_t(t,x) = u_{xx}(t,x)- u(t,x) + g(u(t-h,x)), \ x \in \R, 
\end{equation}
where  $g$ is  a Lipschitzian function such that $|g'|_{L^\infty} >g'(0)$. 
Profile $\phi$ solves the delay
differential equation
\begin{equation}\label{rdl}
 \phi''(t) - c\phi'(t)-\phi(t)+ g(\phi(t-hc))=0, \quad t \in \R.
\end{equation}
Similarly to  Section \ref{nrd} (where we take now $\beta =0$),  we find that  $\phi$ satisfies 
$$
\phi(t)=\mathcal K*g(\phi)(t), \,\,\,
\mathcal K(s)=\frac{1}{\sigma(c)}\left\{\begin{array}{cc} e^{\nu(s-ch)}, & s\geq ch, \\ e^{\mu(s-ch)}, & s<ch,\end{array}\right.
$$
which is exactly the form considered in the DK theory (formally, we set 
 $X=\{\tau\}$, $K(s,\tau)=\mathcal K$ and $g(s,\tau)=g(s)$). Nevertheless, since $L > g'(0)$, the Diekmann-Kaper 
 uniqueness theorem does not apply to (\ref{rdl}).   
 
 In order to use Theorem \ref{mth2}, we first note that  
$$
\chi_1(z,c) =1 - L\int_\R \mathcal K(s)e^{-sz}ds=1-\frac{Le^{-zhc}}{1+cz-z^2}. 
$$
is well defined on $(\nu, \mu)$. Thus,  $\gamma_K=\mu$ and  since $\lim _{x\to \mu-}\int_\R \mathcal K(s)e^{-sx}ds=+\infty$ we obtain that $\gamma_\phi<\gamma_K$. The exponential estimations of $K$ in  {\bf {(SB*)}}, {\bf{ (EC*)}}(ii) are also obviously verified. 

Finally, let $c_\star$ be the minimal value of $c$ for which the equation $z^2-cz-1+Le^{-chz}=0$ has at least one positive root.  This value is well defined and positive. It is easy to see that, for each $c> c_\star$ there exists $m>0$ close to $\lambda_l$ from the right and such that  $\chi_1(m) >0$. 
Hence, we get the following 
  \begin{theorem}\label{lo1}
Suppose that $|g(s)-g(t)|\leq L|t-s|, \, s,t \geq 0,$ and that  $g\in C^{1,\alpha}$ in some neighborhood of $0$ with $g'(0+)>1$. 
Then, for every $c > c_\star$  equation (\ref{rdl}) has at most one  bounded  positive solution $\phi$ vanishing at $-\infty$. 
\end{theorem}
\begin{remark} Theorem \ref{lo1}  gives an alternative proof of  the uniqueness result in  \cite[Theorem 1.1]{atv} where it was additionally assumed that $g\in C^1(\R_+,\R_+)$ and that $g''(0+)$ in finite. Moreover,  we give here a reasonably good  lower bound $c_\star$ for the 'uniqueness' speeds. Observe that if $L =g'(0)$, then   $c_\star$ coincides with the minimal speed of propagation $c_*$.   
\end{remark}

\begin{acknowledgements}
Research was
supported in part by  CONICYT (Chile) through PBCT program ACT-56 and by the University of
Talca, through program ``Reticulados y Ecuaciones". Sergei  Trofimchuk  was partially supported by FONDECYT (Chile), project 1071053. 
\end{acknowledgements}

\end{document}